\documentclass{gtart_h}

\def\ifplaintex{\expandafter\ifx\csname documentclass\endcsname\relax}

\def\gtp{{\mathsurround=0pt\it $\cal G\mskip-2mu$eometry \&\ 
$\cal T\!\!$opology $\cal P\!$ublications}}  

\def\recd{{\small Received:\qua\receiveddate\ifx\reviseddate\relax
\else\qquad Revised:\qua\reviseddate\fi\par}} 

\def\lognumber#1{\def\thelognumber{#1}}
\def\volumenumber#1{\def\thevolumenumber{#1}}
\def\volumeyear#1{\def\thevolumeyear{#1}}
\def\papernumber#1{\def\thepapernumber{#1}}
\def\pagenumbers#1#2{\def\startpage{#1}\def\finishpage{#2}}
\def\published#1{\def\publishdate{#1}}

\def\received#1{\def\receiveddate{#1}}
\def\revised#1{\def\reviseddate{#1}}
\def\accepted#1{\def\accepteddate{#1}}

\long\def\asciiabstract#1{\long\def\theasciiabstract{#1}}

\let\\\par\let\thelognumber\relax\let\thevolumenumber\relax
\let\thepapernumber\relax\let\thevolumeyear\relax\let\startpage\relax
\let\finishpage\relax\let\publishdate\relax\let\receiveddate\relax
\let\reviseddate\relax\let\accepteddate\relax\let\theasciititle\relax
\let\theasciiauthors\relax
\let\theasciiabstract\relax

\let\theasciiemail\relax

\ifplaintex
\font\logobig=cmssbx10 scaled 3836
\font\logomed=cmssbx10 scaled 2557
\else
\font\logobig=cmssbx10 scaled 4200
\font\logomed=cmssbx10 scaled 2800
\fi

\long\def\makeagttitle{   
\count0=\startpage
\agt\hfill      
\hbox to 45truept{\vbox to 0pt{\vglue -13truept{\logomed A\kern -.37em{\logobig 
T}\kern -.38em G}\vss}\hss}
\break
{\small Volume \thevolumenumber\ (\thevolumeyear)
\startpage--\finishpage\nl
Published: \publishdate}

\vglue .25truein

{\parskip=0pt\leftskip 0pt plus
1fil\def\\{\par\smallskip}{\Large\bf\thetitle}\par\medskip} \vglue
0.05truein

{\parskip=0pt\leftskip 0pt plus 1fil\def\\{\par}{\sc\theauthors}
\par\medskip}%
 
\vglue 0.03truein 

{\small\leftskip 25truept\rightskip 25truept{\bf Abstract}\stdspace\theabstract

{\bf AMS Classification}\stdspace\theprimaryclass
\ifx\thesecondaryclass\relax\else; \thesecondaryclass\fi\par
{\bf Keywords}\stdspace \thekeywords\par}\vglue 7truept

}   

\ifplaintex
\font\phead=cmsl9 scaled 950
\font\pnum=cmbx10 scaled 913
\font\pfoot=cmsl9 scaled 950
\headline{\vbox to 0pt{\vskip -4.5mm\line{\small\phead\ifnum
\count0=\startpage ISSN 1472-2739 (on-line) 1472-2747 (printed)
\hfill {\pnum\folio}\else\ifodd\count0\def\\{ }%
\ifx\theshorttitle\relax\thetitle\else\theshorttitle\fi\hfill{\pnum\folio}
\else\def\\{ and }{\pnum\folio}\hfill\ifx\theshortauthors\relax\theauthors
\else\theshortauthors\fi\fi\fi}\vss}}
\footline{\vbox to 0pt{\vglue 0mm\line{\small\pfoot\ifnum\count0=\startpage
\copyright\ \gtp\hfill\else
\agt, Volume \thevolumenumber\ (\thevolumeyear)\hfill\fi}\vss}}
\else
\font=cmsl9 scaled 1050
\font\lnum=cmbx10 
\font=cmsl9 scaled 1050
\makeatletter
\def\@oddhead{{\small\lhead\ifnum\count0=\startpage ISSN 1472-2739 
(on-line) 1472-2747 (printed)\hfill {\lnum\number\count0}\else\ifodd\count0
\def\\{ }\ifx\theshorttitle\relax \thetitle \else\theshorttitle\fi\hfill
{\lnum\number\count0}\else\def\\{ and }{\lnum\number\count0}
\hfill\ifx\theshortauthors\relax 
\theauthors\else\theshortauthors\fi\fi\fi}}\def\@evenhead{\@oddhead}
\def\@oddfoot{\small\lfoot\ifnum\count0=\startpage\copyright\ \gtp\hfill\else
\agt, Volume \thevolumenumber\ (\thevolumeyear)\hfill\fi}
\def\@evenfoot{\@oddfoot}
\makeatother
\fi
\let\maketitlepage\makeagttitle
\let\makeshorttitle\maketitlepage
\let\maketitle\maketitlepage

\newwrite\gtoutfile
\long\gdef\makeheadfile{  
{\def\\{, }\def\s{ }
\immediate\openout\gtoutfile head.xxx
\immediate\write\gtoutfile{Proxy-for: \ifx\theasciiauthors\relax
\theauthors\else\theasciiauthors\fi\s<\ifx\theasciiemail\relax\theemail\else\theasciiemail\fi>}
\immediate\write\gtoutfile{\noexpand\\}
\immediate\write\gtoutfile{Authors: \ifx\theasciiauthors\relax
\theauthors\else\theasciiauthors\fi}
{\def\\{ }\immediate\write\gtoutfile{Title: \ifx\theasciititle\relax
\thetitle\else\theasciititle\fi}}
\immediate\write\gtoutfile{Subj-class: GT or SG, GR etc}
\immediate\write\gtoutfile{MSC-class: \theprimaryclass\ifx\thesecondaryclass\relax\else, \thesecondaryclass\fi}
\immediate\write\gtoutfile{Journal-ref: Algebr. Geom. Topol. \thevolumenumber\s
(\thevolumeyear) \startpage-\finishpage}
\immediate\write\gtoutfile{Comments: Published by Algebraic and
Geometric Topology at}
\immediate\write\gtoutfile{\s\s\s  http://www.maths.warwick.ac.uk/agt/AGTVol\thevolumenumber/agt-\thevolumenumber-\thepapernumber.abs.html}
\immediate\write\gtoutfile{\noexpand\\}
\immediate\write\gtoutfile{}
\ifx\theasciiabstract\relax
\immediate\write\gtoutfile{\theabstract}\else
\immediate\write\gtoutfile{\theasciiabstract}\fi
\immediate\write\gtoutfile{}
\immediate\write\gtoutfile{\noexpand\\}
\immediate\write\gtoutfile{}
\immediate\closeout\gtoutfile}}  

\def\maketitlepage{\makeagttitle\makeheadfile}
\let\makeshorttitle\maketitlepage
\let\maketitle\maketitlepage

\lognumber{13}
\volumenumber{5}
\volumeyear{2005}
\papernumber{13}
\pagenumbers{237}{300}
\received{9 January 2004} 
\revised{16 March 2005}
\accepted{30 March 2005}
\published{15 April 2005}

\usepackage{amsmath,amssymb,amscd,epsf}

\numberwithin{equation}{section} 

\newtheoremstyle{bold}{14pt plus6pt minus6pt}%
{}{\rm}{}{\bf}{}{.75em}%
{\thmname{#1}\thmnumber{#2}\thmnote{\bf\qua #3}}
\theoremstyle{bold}
\newtheorem{sshead}{}
\numberwithin{sshead}{subsection} 
\def\subsubsection#1{\sshead[#1]}
\def\S{section }

\def\cal{\mathcal}

\def\del{\partial}
\def\iso{\cong}
\def\a{\alpha}

\def\D{\Delta}
\def\eps{\epsilon}

\def\s{\sigma}
\def\t{\tau}

\def\nn{\nonumber}

\def\ZZ{\mathbb {Z}/2\mathbb{Z}}

\def\Sn{\mathbb{S}_n}
\def\Sm{\mathbb{S}_m}
\def\Snn{\mathbb{S}_{n+1}}

\def\C{\mathcal{C}}

\def\O{\mathcal{O}}

\def\lab{\mathrm{lab}}

\def\Darc{\mathcal{DA}rc}
\def\Arc{\mathcal{A}rc}
\def\Arcn{\mathcal{A}rc_{\#}}
\def\Tree{\mathcal{T}ree}
\def\Lintree{\mathcal{L}Tree}
\def\Cact{\mathcal{C}act}
\def\Cacti{\mathcal{C}acti}
\def\Corol{\mathcal{S}CC}

\def\Loop{\mathcal{L}oop}
\def\CWCact{K}
\def\UCact{\mathcal{UC}act}

\def\mk{\mathrm{mk}}
\def\wlbptree{\mathcal{T}}
\def\color{cr}

\begin{document}

\title{On several varieties of cacti and their relations}

\author{Ralph M. Kaufmann}
\address
{Department of Mathematics, University of Connecticut\\196
Auditorium Road, Storrs, CT 06269-3009, USA}
\email{kaufmann@math.uconn.edu}

\begin{abstract}
Motivated by string topology and the arc operad, we introduce the
notion of quasi-operads and consider four (quasi)-operads which are
different varieties of the operad of cacti. These are cacti without
local zeros (or spines) and cacti proper as well as both varieties
with fixed constant size one of the constituting loops.  Using the
recognition principle of Fiedorowicz, we prove that spineless cacti
are equivalent as operads to the little discs operad.  It turns out
that in terms of spineless cacti Cohen's Gerstenhaber structure and
Fiedorowicz' braided operad structure are given by the same explicit
chains. We also prove that spineless cacti and cacti are homotopy
equivalent to their normalized versions as quasi-operads by showing
that both types of cacti are semi-direct products of the quasi-operad
of their normalized versions with a re-scaling operad based on
${\mathbb R}_{> 0}$. Furthermore, we introduce the notion of
bi-crossed products of quasi-operads and show that the cacti proper
are a bi-crossed product of the operad of cacti without spines and the
operad based on the monoid given by the circle group $S^1$. We also
prove that this particular bi-crossed operad product is homotopy
equivalent to the semi-direct product of the spineless cacti with the
group $S^1$.  This implies that cacti are equivalent to the framed
little discs operad.  These results lead to new CW models for the
little discs and the framed little discs operad.
\end{abstract}

\asciiabstract{Motivated by string topology and the arc operad, we
introduce the notion of quasi-operads and consider four
(quasi)-operads which are different varieties of the operad of
cacti. These are cacti without local zeros (or spines) and cacti
proper as well as both varieties with fixed constant size one of the
constituting loops.  Using the recognition principle of Fiedorowicz,
we prove that spineless cacti are equivalent as operads to the little
discs operad.  It turns out that in terms of spineless cacti Cohen's
Gerstenhaber structure and Fiedorowicz' braided operad structure are
given by the same explicit chains. We also prove that spineless cacti
and cacti are homotopy equivalent to their normalized versions as
quasi-operads by showing that both types of cacti are semi-direct
products of the quasi-operad of their normalized versions with a
re-scaling operad based on R>0.  Furthermore, we introduce the notion
of bi-crossed products of quasi-operads and show that the cacti proper
are a bi-crossed product of the operad of cacti without spines and the
operad based on the monoid given by the circle group S^1. We also
prove that this particular bi-crossed operad product is homotopy
equivalent to the semi-direct product of the spineless cacti with the
group S^1.  This implies that cacti are equivalent to the framed
little discs operad.  These results lead to new CW models for the
little discs and the framed little discs operad.}



\maketitle

\section*{Introduction}
\addcontentsline{toc}{section}{Introduction}
The cacti operad was introduced by Voronov \cite{V} descriptively as treelike
configurations of circles in the plane to give an operadic
interpretation of the string bracket and Batalin-Vilkovisky (BV)
structure found by Chas and Sullivan \cite{CS} on the loop space
of a compact manifold. The key tool connecting the two is an
``Umkehr'' map in homology by using a Thom-Pontrjagin construction
\cite{V,CJ}. Studying combinatorial models of the moduli space of
bordered surfaces,  we constructed the $\Arc$ operad, which is an
operad built on surfaces with arcs, and showed that this operad
naturally carries compatible structures of a Gerstenhaber (G) and
a BV algebra up to homotopy on the chain level \cite{KLP}.
Moreover the structures above were given by explicit generators
for the operations and explicit homotopies for the relations. The
structure of these generators and relations for the BV operations
bear formal resemblance to those of the cacti operad. With the
help of an additional analysis, we were indeed able to give a map
of the operad of cacti into the $\Arc$ which embeds the former as
a suboperad (up to an overall re-scaling) \cite{KLP} and embeds $\Cacti$ into
the operad $\Darc=\Arc\times \mathbb{R}_{>0}$. This result defined the topology of
$\Cacti$ in terms of metric ribbon graphs. A brief review of this construction
of \cite{KLP} is contained in Appendix B.
 The
operations defining the Gerstenhaber bracket lie in a very
naturally defined smaller suboperad of $\Arc$ than the one
corresponding to cacti. Moreover, adding a generator to this
suboperad which is the element that
 becomes the BV operator we obtain the suboperad corresponding to cacti.
The construction of these spaces and the relation between them is also briefly
reviewed in Appendix B. Going back to string-topology, one
can expect to be able to find a suboperad of
 cacti responsible for the G-bracket. This is the operad of spineless cacti.

Setting these observation in relation to the theorems of F.\ Cohen
\cite{cohen, cohen2} and E.\ Getzler \cite{Ge} ---which state that
Gerstenhaber algebras coincide with algebras over the homology of
the little discs operad and BV algebras coincide with algebras
over the framed little discs operad--- leads us to a comparison of
the operad of cacti and the suboperad mentioned above to the
framed little discs operad and its suboperad of little discs. The equivalence
of spineless cacti and the little discs operad
is one of the main points of this paper. Furthermore, we also give a proof for the Theorem
announced by Voronov \cite{V,Vman}
that cacti   and the framed little discs are equivalent.

One further striking fact about the explicit chain homotopies used
in \cite{KLP} is that they have one more restriction in common,
that of a normalization. This property is, however, not stable
under composition.

The considerations above prompt us to define several different
species of cacti and to study their relations to each other and
their relation to the little discs and the framed little discs
operad. These different species consist of the original cacti,
cacti without additionally marked local zeros which we call cacti
without spines and lastly for both versions their normalized
counterparts, which are made up of circles of radius one. It is
actually a little surprising that one can go through the whole
theory with normalized cacti. For the normalized versions the
gluing rules are slightly different though and are only
associative up to homotopy.

To systematically treat these objects, we introduce the notion of
quasi-operads and define direct, semi-direct and bi-crossed
products of quasi-operads. In this setting, the normalized
versions of cacti and spineless cacti are homotopy associative
quasi-operads, so that their homology quasi-operads are in fact
operads.

In order to define the topological spaces underlying spineless
cacti and cacti, we use a reformulation of the original approach
of \cite{V} in terms of graphs and trees \cite{del}. In this setting, the
spaces for normalized spineless cacti are constructed as
CW-complexes whose cells are indexed by trees. The underlying
spaces for the other versions of cacti are then in turn given as
products of these CW complexes with circles and lines. Moreover
the quasi-operad structure of normalized spineless cacti induced
on the level of cellular chains is already associative and
provides an operad structure \cite{del}.

As to the relation of (spineless) cacti and their normalized versions,
the exact statement is that the non-normalized versions of cacti are
isomorphic as operads to the quasi-operadic semi-direct product of the
normalized version with a scaling operad.
The scaling operad is defined on the
spaces $\Bbb{R}_{>0}^n$ and controls the radii of the circles.
We also show that this semi-direct product is homotopic to the
direct product as quasi-operads.
This makes the normalized and non-normalized versions of (spineless) cacti
 homotopy equivalent as spaces, but
furthermore the products are also compatible up to homotopy, so
that the two versions are equivalent as quasi-operads.

One main result we prove is that the spineless cacti are
equivalent in the sense of \cite{fiedorrecog} to the little discs
operad using the recognition principle of Fiedorowicz
\cite{fiedorrecog}. For the proof, we take up the idea of
\cite{Vman} to use the map contracting the $n+1$-st lobe of a
cactus with $n+1$ lobes. We analyze this map further and prove
that although it is not a fibration that
it is a quasi-fibration. This is done using the Dold-Thom criterium
\cite{DT}.

In this way, the cellular chains of normalized spineless cacti
provide a model of the chains of the little discs operad.
This fact together with indexing of the chains by trees is the
basis for a natural topological solution to Deligne's conjecture
on the Hochschild cohomology of an associative algebra \cite{del}.
Furthermore, in the same spirit the cellular chains of normalized
cacti form an operad which is a model for the framed little discs
which can in turn be used to prove a cyclic version of Deligne's
conjecture. The content of this Theorem \cite{cyclic} is that
there is a cell model of the framed little discs which acts on the
Hochschild complex of a Frobenius algebra.

The theorem about the equivalence of spineless cacti and the
little discs operad and its proof also nicely tie together the
results of \cite{cohen, cohen2} and \cite{fiedorbar,fiedorrecog} in a
geometrical setting on the chain level. By the results of F.\
Cohen \cite{cohen, cohen2} the algebras over the homology of the little
discs operad are Gerstenhaber algebras and by the recognition
principle of Fiedorowicz \cite{fiedorrecog} an operad is
equivalent to the little discs operad if its universal cover is a
contractible braided operad  with a free braid group action. In
our realization, the Gerstenhaber bracket is made explicit on the
chain level. It is in fact given by the signed commutator of a
non-commutative product $*$ which is defined by a path between
point and its image under a transposition under the action of the
symmetric group. This path is also  the path needed to lift the
symmetric group action to a braid action. Moreover, the odd Jacobi
identity of the Gerstenhaber bracket is proved by a relation for
the associator of $*$. In terms of the paths this equation is the
same equation as the braid relation needed to ensure that the
universal cover is a braided operad.

The relationship between framed little discs and little discs is
that the framed little discs are a semi-direct product of the
little discs with the operad built on the circle group $S^1$
\cite{SW}. Actually, for this construction, which we review below,
one only needs a monoid. This example is a special case of a
semi-direct product of quasi-operads.

The relationship for the cacti and spineless cacti is more
involved. To this end, we define the notion of bi-crossed products
of quasi-operads which is an extension of the bi-crossed product
of groups \cite{Ka,Tak}. We show that cacti are a bi-crossed
product of spineless cacti with an operad built on $S^1$. By
analyzing the construction of the operads built on monoids in a
symmetric tensor category where the tensor product is a product,
we can relate the particular bi-crossed structure of cacti to the
semi-direct product with the circle group. To be precise we show
that the particular bi-crossed product giving rise to cacti is
homotopy equivalent to the semi-direct product of cacti without
spines and the group operad built on $S^1$.

The characterization above allows us to give a proof of  the theorem announced by Voronov 
\cite{V,Vman} 
which states that the
operad of cacti is equivalent to the framed little discs operad.
This is done via the equivariant recognition principle \cite{SW}.
Vice-versa the theorem mentioned above together with the characterization of $\Cacti$
as a bi-crossed product imply that
cacti without spines are homotopy equivalent to the little
discs operad.  

Along the way, we give several other pictorial realizations of the
various types of cacti including trees, ribbon graphs and chord
diagrams, which might be useful to relate this theory to other
parts of mathematics. In particular the trees with grafting are
reminiscent of the Connes-Kreimer operads \cite{CK} and in fact as
we have shown in \cite{del} they are intimately related, see also
\S \ref{ckrem} below. The chord diagram approach is close to
Kontsevich's graph realization of the Chern-Simons theory (cf.\
e.g.\ \cite{BN}) and to Goncharov's algebra of chord diagrams
\cite{Go}.

Interpreting the above results inside the $\Arc$ operad, one
obtains that the bi-crossed product corresponding to cacti is
realized as the suboperad corresponding to cacti without spines
and a Fenchel-Nielsen type twist.

The paper is organized as follows:

In the first section, we introduce the notion of quasi-operads
and the operations of forming direct, semi-direct and bi-crossed
products of quasi-operads which we need to describe the cacti
operads and their relations. In \S 2 we then define all the
varieties of cacti we wish to consider, cacti with and without
spines and normalized cacti with and without spines. In addition,
we provide several pictorial realizations of these objects, which
are useful for their study and relate them to other fields of
mathematics. The third section contains the proof that the operad
of spineless cacti is equivalent to the little discs operad. In
paragraph 4, we collect examples and constructions which we
generalize in paragraph 5 in order to study the relations between
the various varieties of cacti.  We start by introducing an operad
called operad of spaces which can be defined in any symmetric
tensor category with products such as topological spaces with
Cartesian product. This operad lends itself to the description of
the semi-direct product with a monoid whose construction we also
review. In the last section, \S5 we then prove that the
non-normalized versions of the (spineless) cacti operads are the
semi-direct products of their normalized version and a
re-scaling operad built on $\Bbb{R}_{>0}^n$. Moreover we show
that this semi-direct product is homotopy equivalent to a direct
product. This section also contains the result that cacti are a
bi-crossed product of cacti without spines and the operad built
on the group $S^1$. Moreover, we show that this bi-crossed
product in turn is homotopic to the semi-direct product of these
operads. These results are then used to give a proof that cacti are equivalent
to the framed little discs.

We also provide two appendices. Appendix A is a compilation of the
relevant notions of graphs and gives the interpretation of cacti
as marked treelike ribbon graphs with a metric. In Appendix B,  we
briefly recall the $\Arc$ operad and the suboperads corresponding
to the various cacti operads and show how to map cacti to elements
of $\Arc$ and vice-versa. This short presentation slightly differs
in style from \cite{KLP}, since we use the language of graphs of
Appendix A to simplify the constructions in the situation at hand.
As such it might also be useful to a reader acquainted with
\cite{KLP}. Furthermore, the $\Arc$ operad provides a
straightforward  generalization of cacti to higher genus and even
allows to additionally introduce punctures.

\subsection*{Acknowledgments}
\addcontentsline{toc}{subsection}{Acknowledgments}

It is a pleasure to thank A.A.\ Voronov for
discussions. We are also indebted to him for sharing the
manuscript \cite{Vman}. Our special thanks also goes to J.\ E.\
McClure for discussions about the recognition principle of
Fiedorowicz and to the referee for pointing out needed
improvements. We especially  thank  the Max-Planck-Institut f\"ur
Mathematik in Bonn where a great part of the work for this paper
was carried out for its hospitality and support. This work
received partial support of the NSF under grant \#0070681 which is
also gratefully acknowledged.

\section{Quasi-operads and direct, semi-direct and\newline bi-crossed products}
\label{quasiopsection}

In our analysis of the various types of  cacti, we will need a
structure which is slightly more relaxed than operads. In fact, in
the normalized versions of cacti the compositions will fail to be
associative on the nose, although they are associative up to
homotopy. This leads us to define and study quasi-operads. These
quasi-operads afford  certain constructions such as semi-direct
products and bi-crossed products which are not necessarily defined
for operads. On the other hand semi-direct products and bi-crossed
products of quasi-operads may yield operads. If one is mainly
interested in the homology operads, it is natural to consider
quasi-operads which are associative up to homotopy. Lastly, in
certain cases quasi-operads can already provide operads on the
chain level as our normalized spineless cacti below \cite{del}.

\subsection{Quasi-operads}

We fix a strict monoidal category $\mathcal{C}$ and denote by
$\Sn$ the symmetric group on $n$ letters. A  quasi-operad is an
operad where the associativity need not hold. More precisely:

\subsubsection{Definition}
\label{quasiopdef}
A quasi-operad $\mathcal C$ is a collection of objects
$\O:=\{O(n):O(n) \in \mathcal{C}, n\geq1\}$ together with an $\Sn$
action on $O(n)$ and maps called compositions
\begin{equation}
\circ_i: O(m) \otimes O(n) \rightarrow O(m+n-1), i \in \{1,\dots
m\}
\end{equation}
 which are  $\Sn$-equivariant: if $op_m \in O(m)$ and $op_n \in O(n)$
\begin{equation}
\label{snequivar}
\s_m(op_m) \circ_i \s'_n(op_n) =
\s_m\circ_i\s'_n(op_m\circ_{\s_m(i)} op_n)
\end{equation}
where $\s_m\circ_i\s'_n\in \mathbb{S}_{m+n-1}$ is the permutation
that the block or iterated permutation
\begin{multline}
(1,2,\dots,i-1,(1',\dots,m'), i+1 \dots, n) \mapsto\\
\s_n(1,2,\dots,i-1,\s'_m(1',\dots,m'), i+1 \dots, n)
\end{multline}
induces on $(1'',\dots,(m+n-1)'')$ where
$$
j'' = \begin{cases} j& 1\leq j \leq i-1\\
 j-i+1'&i\leq j \leq  i+n-1\\
j-n&i+n \leq j \leq m+n-1
\end{cases}
$$
A quasi-operad is called unital if an element $id \in O(1)$ exists
which satisfies for all $op_n \in O(n), i\in \{1,\dots,n\}$
$$
\circ_i(op_n,id)= \circ_1(id,op_n)= op_n
$$

\subsubsection{Remark}
If a  quasi-operad  in the topological category is homotopy
associative then its homology has the structure of an operad.
 In certain cases, like the ones we will consider, the structure of
 an operad
already exists on the level of a chain model.

\subsubsection{Definition} A morphism of quasi-operads is a map which
preserves all structures.

An isomorphism of quasi-operads is an invertible morphism of
quasi-operads. This will be denoted by $\iso$.

A quasi-operad morphism $A \rightarrow B$ is said to be an
equivalence if for each $k \geq 0$, $A(k) \rightarrow B(k)$ is a
$\mathbb{S}_k$-equivariant homotopy equivalence. This relation
will be denoted by $\simeq$.

\subsubsection{Definition}
Two quasi-operad structures $\circ_i$ and $\circ'_i$ on a fixed
collection of $\Sn$-spaces $O(m)$  are called homotopy equivalent
(through quasi-operads) denoted by $\sim$ if there is a homotopy
of  maps
\begin{equation}
\circ_i(t): O(m) \otimes O(n) \rightarrow O(m+n-1)
\end{equation}
for $t\in [0,1]$ such that \ $\circ_i(0)=\circ_i,
\circ_i(1)=\circ'_i$ and such that for any fixed $t$ the
$\circ_i(t)$ give the $O(n)$ the structure of a quasi-operad.

\subsubsection{Remark} Two homotopy equivalent quasi-operads induce
isomorphic structures on the homology level.

\subsubsection{Definition} An operad is a quasi-operad for which
associativity holds, i.e.\ for $op_k \in O(k), op'_l \in O(l)$ and
$op''_m \in O(m)$
\begin{equation}
\label{assequation}
(op_k \circ_i op'_l) \circ_j op''_m=
\begin{cases}
(op_k \circ_j op''_m) \circ_{i+m-1} op'_l&\text{ if } 1\leq j < i \\
op_k \circ_i (op'_l \circ_{j-i+1} op''_m)&\text{ if } i\leq j< i+l\\
(op_k \circ_{i-l+1} op'_l) \circ_j op''_m&\text{ if } i+l\leq j
\end{cases}
\end{equation}

An operad morphism is a map of collections preserving all the
operad structures.

\subsubsection{Remark}
Note that our operads correspond to the pseudo-operads of
\cite{MSS}. In case a unit exists these two notions coincide
\cite{MSS}. We drop the ``pseudo'' in our nomenclature in order to
avoid confusion between quasi- and pseudo-operads. In a strict
sense, our quasi-operads are quasi-pseudo-operads which is
certainly an expression we wish to avoid.

\medskip

 We will use the following terminology of
\cite{fiedorrecog}.

\subsubsection{Definition}
\label{enterm}
 An operad morphism $A \rightarrow B$ is said to be an
equivalence if for each $k \geq 0$, $A(k) \rightarrow B(k)$ is a
$\mathbb{S}_k$-equivariant homotopy equivalence. We say that an
operad $A$
 is $E_n$
($n = 1, 2, 3, \dots, \infty$) if there is a chain of operad
equivalences (in either or both directions) connecting $A$ to the
Boardman-Vogt little $n$-cubes operad $C_n$ (cf. \cite{BV}).

\subsection{Direct products}
\subsubsection{Definition} Given two quasi-operads
$C(n)$ and $D(n)$ in the same category, we define their direct
product $C\times D$ to be given by  $(C\times D)(n):= C(n)\times
D(n)$ with the diagonal $\Sn$ action, i.e.\ the action of $\Sn$
induced by the diagonal map $\Sn \rightarrow \Sn \times \Sn$, and
the compositions
$$
(c,d) \circ_{i,\C\times D} (c',d'):=(c\circ_{i,C} c',d\circ_{i,D}
d')
$$
Since the compositions are componentwise it follows that:
\subsubsection{Proposition}
{\sl The direct product of two operads is an operad.}

\subsection{Semi-direct products}

\subsubsection{Definition}
Fix two  quasi-operads $C(n)$ and $D(n)$ in the same category
together with a collection of morphisms:
\begin{eqnarray}
\label{semiops} \circ_i^D: C(n)\times D(n) \times C(m)
&\rightarrow& C(n+m-1)
\quad i = 1, \dots n\nn\\
(c,d,c')&\mapsto&\circ_i^D(c,d,c')=:c \circ_i^d c'
\end{eqnarray}
which satisfy the analog of equation (\ref{snequivar}) for the action of $\Sn
\times \Sm$ with $\Sn$ acting diagonally on the first two factors.

Note that we used the superscripts to indicate that we view the
dependence of the map on $D$ as a perturbation of the original
quasi-operad structure.

We define  the semi-direct product $C\rtimes D$ with respect to
the $\circ_i^D$ to be given by the collection $(C \times D) (n):=
C(n)\times D(n)$ with diagonal $\Sn$ action and compositions
\begin{eqnarray}
\label{semiprod}
 \circ_i: (C\times D)(n) \times (C \times D)(m) &=&
C(n)\times D(n) \times C(m) \times D(m)\nn\\
\longrightarrow&& C(n+m-1)\nn\\
(c,d) \circ_i(c',d')&=& (c\circ_i^d c',d \circ_i d')
\end{eqnarray}
where we use the upper index on the operations to show that we use
the universal maps (\ref{semiops}) with fixed middle argument.

In  case we are dealing with unital (quasi) operads we will also
require  that $1_C \times 1_D$ is a unit in the obvious notation,
and that $\circ_i^{1_D}=\circ_i$.

\subsubsection{Remark} In general the semi-direct product of two operads
need not be an operad. This depends on the choice of the
$\circ_i^D$. Of course the direct product of two operads is a
semi-direct product of quasi-operads which is an operad. If the
quasi-operad $D$ fails to be associative, then any semi-direct
product $C \rtimes D$ will fail to be an operad. However if $D$ is
an operad and $C$ is just a quasi-operad, it is possible that
there are maps $\circ_i^D$ s.t.\ $C \rtimes D$ will be an operad.

\subsubsection{Definition} We call a quasi-operad $C$  {\em
normal} with respect to an operad $D$ and maps of the type
(\ref{semiops}) if the semi-direct product $C\rtimes D$ with
respect to these maps yields an operad, i.e.\ is associative.

\subsubsection{Examples}
Examples of this structure are usually derived if the operad $D$
acts on the quasi-operad $C$ and the twisted operations are
defined by first applying this action. The semi-direct product
of an operad with a monoid \cite{SW} is such an example (see below)
as are the semi-direct products
of \cite{KLP} and the (spineless) cacti with respect to their
normalized versions (see below).

A sometimes useful criterion is:
\subsubsection{Lemma}
\label{actionlemma} {\sl Consider two operads $C$ and $D$ with  left
actions of $D(n)$ on $C(m)$,
$$
\rho_i: D(n) \times C(m) \rightarrow C(m) \quad i \in \{1,\dots, n\}
$$
s.t.\ the maps
\begin{eqnarray*}
\circ_i^D: C(n)\times  D(n) \times C(m) &\rightarrow& C(m):\\
(c,d,c') &\mapsto& c \circ_i \rho_i(d)c'
\end{eqnarray*}
are $\Sn \times \Sm$ equivariant in the sense of \ref{quasiopdef}
with $\Sn$ acting diagonally on the first two factors and
for $c' \in C(l), i\leq j< i+l$
\begin{equation}
\label{semiopcond}
 \rho_i(d)(c') \circ_{j-i+1} (\rho_j(d \circ_i d') c'') =
 \rho_i(d)[c'\circ_{j-i+1} \rho_{j-i-1}(c'')]
\end{equation}
 then $C \rtimes D$ with
respect to $c \circ_i^d c'$ is an operad.}

\begin{proof}
Since the action is compatible with the $\Sn$ actions it remains
to check the associativity. The first and third case of the
equation (\ref{assequation}) clearly hold and the second case
follows from the equation (\ref{semiopcond}).
\end{proof}

\subsubsection{The right semi-direct product}
\label{rsemi} There is a right version of the semi-direct
products using maps
\begin{eqnarray}
\label{rsemiops} \circ^C_i: D(n)\times C(m) \times D(m)
&\rightarrow& D(n+m-1)
\quad i = 1, \dots n\nn\\
(d,c,d')&\mapsto&\circ_i^C(d,c,d')=:d\circ_i^c d'
\end{eqnarray}
which again define compositions on the products $(C\times D) (n):=
C(n)\times D(n)$ via
\begin{eqnarray}
\label{rsemiprod} \circ_i: (C\times D)(n) \times (C \times D)(m)
&=&
C(n)\times D(n) \times C(m) \times D(m)\nn\\
\longrightarrow&& C(n+m-1)\nn\\
(c,d) \circ_i(c',d')&=& (c\circ_i  c',d \circ_i^{c'} d')
\end{eqnarray}
We call a quasi-operad $D$  {\em normal} with respect to $C$ if
there are maps of the type (\ref{rsemiops}) such that the maps
(\ref{rsemiprod}) give an operad structure to the product $C
\times D$. Given an operad $C$ an a quasi-operad $D$ normal with
respect to $C$, we call the product together with the operad
structure (\ref{semiprod}) the {\em right semi-direct product} of
$C$ and $D$ which we denote $C\ltimes D$.

This structure typically appears when one has a right action of
the operad $C$ on the quasi-operad $D$.

\subsection{Bi-crossed products}
\label{bicrossedsection}

\subsubsection{Definition}
Consider two quasi-operads $C(n)$ and $D(n)$ together with a
collection of maps:
\begin{eqnarray}
\label{biops} \circ_i^D: C(n)\times D(n) \times C(m) &\rightarrow&
C(n+m-1)
\quad i = 1, \dots n\nn\\
(c,d,c')&\mapsto&\circ_i^D(c,d,c')=:c \circ_i^d c'\\
\label{biops2} \circ_i^C: D(n) \times C(m) \times
D(m)&\rightarrow& C(n+m-1) \quad i = 1, \dots n\nn\\
(d,c,d')&\mapsto&\circ_i^C(d,c,d')=:d \circ_i^c d'
\end{eqnarray}
where the operations (\ref{biops}) satisfy the analog of equation
(\ref{snequivar}) for the action of $\Sn \times \Sm$ with $\Sn$
acting diagonally on the first two factors and the operations
(\ref{biops2}) satisfy the analog of equation (\ref{snequivar})
for the action of $\Sn \times \Sm$ with $\Sm$ acting diagonally on
the second two factors.

Again we used the superscripts to indicate that we view the
dependence on the other quasi-operad as a perturbation of the
original quasi-operad structure.

We define the bi-crossed product $C\bowtie D$ with respect to the
operations $\circ_i^D, \circ_j^C$ to be given by the collection $(C
\times D) (n):= C(n)\times D(n)$ with diagonal $\Sn$ action and
compositions
\begin{eqnarray}
\label{biprod}
 \circ_i (C\times D)(n) \times (C \times D)(m) &=&
C(n)\times D(n) \times C(m) \times D(m)\nn\\
\longrightarrow&& C(n+m-1)\nn\\
(c,d) \circ_i(c',d')&=& (c\circ_i^d c',d \circ_i^{c'}d')
\end{eqnarray}
where we use the upper index on the operations to show that we use
the universal maps (\ref{biops}) with fixed middle argument.

In the case we are dealing with unital operads we will also
require that the perturbed compositions are such that $1_C \times
1_D$ is a unit in the obvious notation, and that
$\circ_i^{1_C}=\circ_i$ and $\circ_i^{1_D}=\circ_i$ .

\subsubsection{Remark} Again it depends on the choice of the
$\circ_i^D, \circ_j^C$ if the bi-crossed product of two operads
is an operad. In the case that $\circ_i^D, \circ_j^C$ are given by
actions as in Lemma \ref{actionlemma} then this is guaranteed if
the condition (\ref{semiopcond}) and its right analog hold.

Furthermore it is possible that the bi-crossed product of two
quasi-operads is an operad.

\subsubsection{Definition} We call two quasi-operads $C$ and $D$ {\em
matched} with respect to maps of the type (\ref{biops}) and
(\ref{biops2}) if the quasi-operad  $C \bowtie D$ is an operad.

\subsubsection{Examples}$\phantom{99}$
\begin{itemize}
    \item [1)] Factoring the multiplication maps through the first
    and third projection and using the structure maps of the two operads
    we obtain the direct product.
    \item [2)] If $D$ is the operad based on a monoid and choosing
    the maps $\circ_i^C$ to be unperturbed and defining the maps
    $\circ_i^D$ as in (\ref{semi}) we obtain the semi-direct product.
    \item [3)] If $C$ and $D$ are concentrated in degree $1$ and happen
    to be groups then the bi-crossed product is that of matched
    groups \cite{Ka,Tak}. Otherwise we obtain that of matched monoids.
    \item[4)] Below we will show that the cacti operad is a
    bi-crossed product of the cactus operad and the operad built
    on $S^1$ as a monoid and thus these operads are matched with respect
    to the specific perturbed multiplications given below.
\end{itemize}

\section{Several varieties of cacti}
\subsection{Introduction}
The operad of $\Cacti$ was first introduced by Voronov in \cite{V}
as pointed treelike configurations of circles. In the following we
will first take up this description and introduce Cacti and
several versions of related operads in this fashion. This approach
is historical and lends itself to describe actions on the loop
space of a compact manifold \cite{CJ,V} of cacti. For other
purposes, especially giving the topology on the space of cacti,
other descriptions of a more combinatorial nature are more
convenient. In the following section, we will both recall the
traditional approach as well the sometimes more practical
definition in terms of graphs.

\label{defs}
 \nopagebreak
\subsection{General setup for configurations of circles}
\nopagebreak
\subsubsection{Notation} By an $S^1$ in the plane we will mean
a map of the standard $S^1 \subset
\Bbb{R}^2$ with the induced metric and orientation which is  an
orientation preserving  embedding $f:S^1 \rightarrow \Bbb
{R}^2$. A configuration of $S^1$'s in the plane is given by a collection of
 finitely many of
these maps which have at most finitely many  intersection points
in the image. I.e.\  if $f_i:1\leq i \leq n$ is such a collection,
then if $i\neq j:$ $f_i(\theta_i)=f_j(\theta_j)$ for only finitely
many points $\theta_i$.

By an $S_r^1$ in the plane we will mean a map of the standard
circle of radius $r$: $S_r^1 \subset \Bbb{R}^2$ with the induced
metric and orientation which is an orientation preserving
embedding $f:S_r^1 \rightarrow \Bbb {R}^2$. A configuration of
$S^1_r$s in the plane is a collection of finitely many of these
maps which have at most finitely many  intersection points in the
image.

\subsubsection{Re-parameterizations}
Notice that a circle in a plane comes with a natural parameter.
This is inherited from the natural parameter $\theta_r$ of the
standard parametrization of $S^1_r:
(r\cos(\theta_r),r\sin(\theta_r))$. Sometimes we have to
re-parameterize a circle, so that its length changes. To be
precise let $f:S^1_r \rightarrow \Bbb{R}^2$ be a parametrization
of a circle in the plane. Then $f_R: S^1_R \rightarrow \Bbb{R}^2$,
called the re-parametrization to length $R$, is defined to be the
map $f_R= f\circ rep_r^R$ with $rep_r^R:S^1_R \rightarrow S^1_r$
given by $\theta_R \mapsto \theta_r$.

\subsubsection{Dual black and white graph}
Given a configuration of $S^1$s in the plane, we can associate to
it a dual graph in the plane. This is a graph with two types of
vertices, white and black. The first set of vertices is given by
replacing each circle by a white vertex. The second set of
vertices is given by replacing  the intersection points with black
vertices. The edges run only from white to black vertices, where
we join two such vertices if the intersection point corresponding
to the black vertex lies on the circle represented by the white
vertex.

We remark that all the $S^1$s are pointed by the image of $0=(0,1)$.
On any given $S^1$ in the plane we will call this point and its image
local zero or base point.

If we are dealing with circles of radii different from one we will
label the vertices of the trees by the radius of the respective
circle.

\subsubsection{Cacti and trees} The configurations corresponding to
cacti will all have trees as their dual graphs. We would like to
point out that these trees are planar trees, i.e. they are
realized in the plane or equivalently they have a cyclic order of
each of the sets of edges emanating from a fixed vertex.

Moreover we would like to consider rooted trees, i.e.\ trees with
one marked vertex called root. Recall that specifying a root
induces a natural orientation for the tree and a height function
on vertices. The orientation of edges is toward the root and the
height of a vertex is the number of edges traversed by the unique
shortest path from the vertex to the root. Due to the orientation
we can speak of incoming and outgoing edges, where the outgoing
edge is unique and points toward the root. Naturally the edges
point from the higher vertices to the lower vertices. A rooted
planar tree has a cyclic order for all edges adjacent to a given
vertex and a linear order on the adjacent edges to any given
vertex except the root. A planar rooted tree with a linear order
of the incoming edges at the root is called a planted tree. The
leaves of a rooted tree  are the vertices which only have outgoing
edges.

We call a configuration of $S^1$s in the plane rooted if one of
the circles is marked by a point, we also call this point the
global zero. In this case, we include a black vertex in  the dual
graph for this marked point and make this the root, so that the
dual graph of a rooted configuration is a rooted tree. This tree
is actually also planted by the linear order of the incoming edges
of the root provided by making the component on  which the root
lies the smallest element in the linear order.

Given such a rooted configuration of $S_r^1$s in the plane we call
the images of the $0$s ($f_i(0)$) together with the image of the
marked point and the intersection points the special points of the
configuration. We also call the connected components of the image
minus the special points the arcs of the configuration. If
$f_i|(\theta_1, \theta_2)=a$ then we define
$|a|:=\frac{1}{2\pi}(\theta_2-\theta_1)$ to be the length of $a$.
In the same situation we define $\bar a:=f_i|[\theta_1, \theta_2]$
to be the closure of $a$.

\subsubsection{Definition}
Given a  configuration of $S^1$s in the plane whose dual graph is
a  tree, we say that an $S^1$ is contained in another $S^1$ if the
image first circle is a subset of the disc bounded by the image of
the second circle.

A configuration of $S^1$s in the plane is called {\it tree-like}
if its dual graph is a connected tree and no $S^1$ is contained in
another $S^1$.

\subsubsection{The perimeter or outside circle and the global zero}
\label{perimeter}
For a marked tree-like configuration of $S_r^1$s in the plane let
$R=\sum r_i$ then there is a surjective map of $S_R^1$ to the
image of this configuration which is a local embedding. This map
is given by starting at the marked point or global zero of the
root of the configuration going around this circle in the positive
sense until one hits the first intersection point and then
starting to go around the next circle in the cyclic order again in
the positive direction until the next intersection point and so on
until one again reaches the zero of the root.

We will call this map and, by abuse of notation, its image the
perimeter or the outside circle of the configuration. We will also
call the zero of the perimeter the global zero.

\subsection{Normalized cacti without spines}

We will now introduce normalized spineless cacti as configurations. Later,
to give a topology, we will reinterpret these configurations in terms of graphs.

\subsubsection{Definition}
We define {\em normalized cacti without spines}  to be labelled rooted
collections of parameterized $S^1$s grafted together in a tree like
fashion with the gluing points being the zeros of the $S^1$. More
precisely we set:

$\Cact^1(n) :=\{$rooted tree-like configurations of $n$ labelled
$S^1$s in the plane such that the root (global zero) coincides
with the marked point (zero) of the component it lies on,  the
points of intersection are such that the circles of greater height
all have zero as their point of intersection.$\}$/ isotopies
preserving the incidence conditions.

Here and below preserving the incidence conditions means that if
$f_{i,t}:S^1_{r_i}\times I$ are the isotopies and
$f_{i,0}(p)=f_i(p)=f_j(q)=f_{j,0}(q)$ then for all $t$:
$f_{i,t}(p)=f_{j,t}(q)$ and vice-versa if $f_{i,0}(p)=f_i(p)\neq
f_j(q)=f_{j,0}(q)$ then for all $t$: $f_{i,t}(p)\neq f_{j,t}(q)$.

\medskip

We will take the conditions ``without spines'' and ``spineless''
to be synonymous.

\subsubsection{Remark}
We would like to point out that there is only one zero which is
not necessarily an intersection point, namely that of the root. It
can however also be an intersection point. Moreover, one could
rewrite the condition of having a dual black and white graph that
is a tree in the form: given any two circles their intersection is
at most one point.

\subsubsection{Remark} There are several ways to give a topology to
this space. One way to give it a topology is by describing the
degenerations of the above configurations, as was done originally
in \cite{V}. This is done by allowing the intersection points and
the root to move in such a way that they may collide, and ``pass''
each other moving along on the outside circle. If an intersection
point collides with the marked point from the positive direction,
i.e.\ the length of the arc going counterclockwise from the root
to the intersection point goes to zero, then the root passes to
the new component.

The quickest way is to give the topology to the spaces $\Cact^1$
as subspaces of the operad $\Darc$ as defined in \cite{KLP}. A
brief review of the necessary constructions is given in the
Appendix B below for the reader's convenience.

 Lastly, one can define the topology in terms of combinatorial data
 by gluing of products of simplices indexed by trees
 as first explained in \cite{del}.
For definiteness, we will use this construction to fix our
definitions.

\subsubsection{Notation}
Recall that the dual tree of a cactus is a bi-colored (b/w)
bi-partite planar planted tree. Such a tree has a natural
orientation towards the root. We call an edge white if it points
from a black to a white vertex in this orientation and call the
set of these white edges $E_w$. Let  $V_w$ be the set of white
vertices. For a vertex $v$ we let $|v|$ be the set of incoming
edges, which is equal to the number of white edges incident to $v$
and is also equal to the total number of edges incident to $v$
minus one. We call $\wlbptree(n)$ the set of planar planted
bi-partite trees with white leaves, black root, and $n$ white
vertices which are labelled from $1$ to $n$.

\subsubsection{Definitions} The topological type of a spineless
normalized cactus in $\Cact^1(n)$ is defined to be the tree $\t
\in \wlbptree(n)$ which is its dual b/w planar planted tree
together with the labelling of the white vertices induced from the
labels of the cactus.

We define $\wlbptree(n)^k$ to be the elements of $\wlbptree(n)$
with $|E_w|=k$.

Let $\Delta^n$ denote the standard $n$-simplex, $ |\D^n|$ its
standard realization in $\mathbb{R}^{n+1}$ as $\{(t_1,\dots t_{n+1})|\sum_i t_i=1\}$.
We denote the interior
of $|\D^n|$ by $|\dot \D^n|$.

For $\t \in \wlbptree$ we define

\begin{equation}
\label{deltataudef}
\D(\t):=\times_{v \in V_w(\tau)}|\D^{|v|}|
\end{equation}
Notice that $\dim(\D(\t))=|E_w(\t)|$ and that the set $E_w$ has
a linear order which defines an orientation of $\D(\t)$.

We also let
\begin{equation}
\dot \D(\t):=\times_{v \in V_w(\tau)}|\dot \D^{|v|}|
\end{equation}

\subsubsection{Lemma}{\sl
A normalized spineless cactus is uniquely determined by its
topological type and the length of the arcs.}

\begin{proof}
It is clear that each normalized spineless cactus gives rise to
the described data. Vice versa given the data, one can readily
construct a  representative of a spineless cactus with the
underlying data. A quick recipe is as follows. Realize the given
tree in the plane. Blow up the white vertices to circles which do
not intersect and do not contain any of the black points.  Mark
the segments of the circles between the edges by the label of the
second edge bounding the arc in the counterclockwise orientation.
Mark the intersection point of the first edge in the linear order
of the root with the unique circle it intersects. Delete the part
of the edges inside these circles. Now contract the edges and if
necessary deform the circles during the contraction such that they
do not touch. This is only a finite problem and thus such choices
can be made. We will call the images of the circles lobes. The
lobes are labelled from 1 to n by the label of the vertices. There
are obvious maps of $S^1$ onto each lobe, which have lengths of
arcs between the special points (intersection or marked)
corresponding to the labelling. This gives a representative. Since
the data is invariant under isotopy preserving the intersections,
we have constructed the desired cactus and hence shown the
bijection.
\end{proof}

\subsubsection{Lemma}{\sl
For a normalized spineless cactus the lengths of the segments
lying on a given lobe represented by a vertex $v$ are in 1-1
correspondence with points of the open simplex $\dot
\Delta^{|v|}$.}

\begin{proof}
The lengths of the arcs have to sum up to the radius of the lobe
which is one and the number of arcs on a given lobe is $|v|+1$.
\end{proof}

\subsubsection{Remark} The Lemma above also gives an identification
of the arcs of a cactus with the edges of the tree $\t$ specifying
its topological type. Here we fix that an edge $e$ incident to a
white vertex $v_w$ and a black vertex $v_b$ corresponds to the arc
on the lobe of $v_w$ running from the special point (intersection
or root) preceding $v_b$ to $v_b$.

\subsubsection{Proposition}
\label{CWcactset}
{\sl As sets $\Cact^1(n)=\amalg_{\t\in
\wlbptree(n)} \dot\D(\t)$.}
\begin{proof}
Immediate by the preceding two Lemmas.
\end{proof}

\subsubsection{Degenerations}
By the above we can identify the vertices of $\D(\t)$ and
therefore the coordinates of points of $\D(\t)$ with the arcs of a
cactus of topological type $\t$.

Given a cactus $c$ and an arc $a$ with length $|a|<1$ of $c$ we
define the degeneration of $c$ with respect to $a$ to be the
configuration of $S^1$s obtained by a homotopy contracting the
closure of the arc $\bar a$ to a point $p$, but preserving all
other incidence conditions, together with the following root. If
$a$ is not the first arc on the outside circle, then the root
remains unchanged. If $a$ is the first arc of the outside circle
and $a'$ is the second arc of the outside circle which lies in the
image of the map $f_j$, then the new root is defined to be the
point which is the pre-image of $p$ under the map $f_j$, i.e.\
$f_j^{-1}(p)$.

\subsubsection{Degeneration of trees} There is also a purely
combinatorial way to describe the degeneration of the b/w
bi-partite planar planted tree by cutting and re-grafting
\cite{del}. An abbreviated non-technical version is as follows.
Given $\t \in \wlbptree(n)$ and an edge $e$ in $\t$ incident to a
white vertex $v_w$ with $|v_w|>0$ the contraction of $\t$ with
respect to $e$ is given by the following procedure. First let
$v_w,v_b$ be the white and black vertices $e$ is incident to and
let $e'$ be the edge immediately preceding $e$ in the cyclic order
at $v_w$. Let $v_w$ and $v_b'$ be the vertices of $e'$. The
contraction of $\t$ with respect to $e$ is given by the tree in
which the edge $e$ and the vertex $v_b$ are removed and the
remaining branches of $v_b$ are grafted to $v'_b$ in such a way
that they keep their linear order and immediately precede the edge
$e'$ in the cyclic order at $v'_b$.

\subsubsection{Remark}
The reader can readily verify that the degeneration of cacti in
the arc and combinatorial interpretations agree.

\subsubsection{A CW-Complex}
Given a cell $\D(\t)$ and a vertex $v$ of any of the constituting
simplices of $\D(\t)$ we define the $v$-th face of $\D(\t)$ to be the subset of
$\D(\t)$ whose points have $v$-th coordinate equal to zero.

We let $\CWCact(n)$ be the CW complex whose k-cells are indexed by
$\t \in \wlbptree(n)^k$ with the cell $C(\t)=|\D(\t)|$ and the
attaching maps $e_{\t}$ defined as follows. We identify the $v$-th
face of $\D(\t)$ with $\D(\t')$ where $\t'$ is the topological
type of the cactus $c'$ which is the degeneration of a cactus $c$
of topological type $\t$ with respect to the arc $a$ that
represents the vertex $v$.

We denote by $\dot e_{\t}$ the restriction of $e_{\t}$ to the
interior of $\D(\t)$. Notice that $\dot e_{\tau}$ is a bijection.

\subsubsection{Theorem}
\label{cact1} {\sl The elements of $\Cact^1(n)$ are in bijection with
the elements of the CW complex $\CWCact(n)$.}

\begin{proof}
Immediate from the Proposition \ref{CWcactset} above.
\end{proof}

\subsubsection{Definition}
 We will use the above theorem to give $\Cact^1$ the
topology induced by the above bijection, that is we define the
topological space $\Cact^1(n)$ as
$$\Cact^1(n):=\CWCact(n).$$

\subsubsection{The action of  $\Sn$}
There is an action of $\Sn$ on $\Cact^1(n)$ which acts by
permuting the labels.

\subsubsection{Gluing}
\label{glue}
We define the following operations
\begin{equation}
\circ_i: \Cact^1(n) \times \Cact^1(m) \rightarrow \Cact^1(n+m-1)
\end{equation}
by the following procedure: given two normalized cacti without
spines we re-parameterize the $i$-th component circle of the first
cactus to have length $m$ and glue in the second cactus by
identifying the outside circle of the second cactus with the
$i$-th circle of the first cactus.

These gluings do not endow the normalized spineless cacti with the
structure of an operad, but with the slightly weaker structure of
a quasi-operad of section \ref{quasiopsection}.

By straightforward computation we have the following:
\subsubsection{Proposition}
{\sl The glueings make the spaces $\Cact^1(n)$ into a topological
quasi-operad.}

\subsubsection{Remark} 
\begin{figure}[ht!]
\epsfxsize = \textwidth
\epsfbox{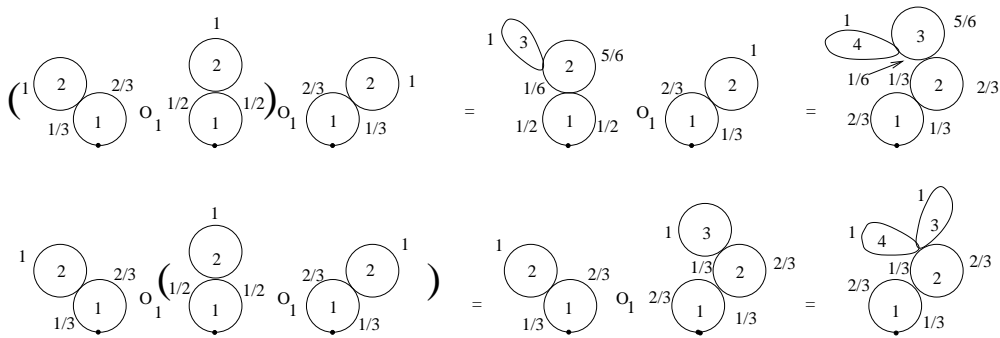}
\caption{\label{quasicact}
An example for non-associativity in $\Cact^1$}
\end{figure}

The above gluing operations are indeed not strictly associative
as the example in Figure \ref{quasicact}
 shows. As in this example, the gluings are
associative up to homotopy  in general, as we will discuss below.

\subsubsection{Remark}
It is shown in \cite{del} that there is an operad structure on
normalized spineless cacti on the cellular chains of $\CWCact$.

This together with the Theorem \ref{cactdisc} provides the basis for a
 proof of Deligne's conjecture on the Hochschild
cohomology of an associative algebra \cite{del}.

\subsection{Cacti without spines}

\subsubsection{Definition}
We define {\em cacti without spines} by an analogous procedure to
that of normalized cacti only this time taking $S^1_r$s, i.e.\
circles of different radii.

As a set $ \Cact(n) :=\{$rooted tree-like configurations of $n$
labelled $S_r^1$s in the plane such that the root (global zero)
coincides with the marked point (zero) of the component it lies on
and that the points of intersection are such that the circles of
greater height all have zero as their point of intersection $\}$/
isotopies preserving the incidence conditions.

Again $\Sn$ acts via permuting the labels.

\subsubsection{Lemma}
$\Cact(n)=\Cact(n)^1\times \mathbb{R}^n_{>0}$.

\begin{proof}
As in the previous case of normalized spineless cacti, such a
configuration is given bijectively by its topological type and the
lengths of its arcs. Each arc belongs to a unique lobe and the sum
of the lengths of the arcs belonging to a lobe is the radius of
the given $S^1$. Let $\mathbf{l}_i:=(l_{i_1},\dots ,l_{i_s})$ be the collections of
lengths of the arcs of the lobe $i$ whose radius is $r_i=\sum_j l_{i_j}$ then
$\mathbf{l}_i$ corresponds to a unique point  in
$\Delta^{s-1}\times \mathbb{R}_{>0}$ given by
$((l_{i_1}/r_i,\dots,l_{i_s}/r_i),r_i)$. This establishes the
claimed bijection.
\end{proof}

\subsubsection{Definition}
As a topological space, we define $\Cact(n)$ to be
$$\Cact(n):=\Cact(n)^1\times
\mathbb{R}^n_{>0}$$
with the product topology.

\subsubsection{Remark}
A description of the topology on this space is given, by
allowing the intersection points and the global zero to move and
collide and pass each other along the outside circle as before,
with the same rule for the global zero as before and also letting
the radii vary. This topology agrees with the
product topology $\Cact(n)=\Cact^1(n)\times \mathbb{R}_{>0}^{n}$ above.

It also agrees with the one induced by the embedding of the
operad into $\Darc$ \cite{KLP}, see also Appendix B.

\subsubsection{Gluing} We define the following operations
\begin{equation}
\circ_i: \Cact(n) \times \Cact(m) \rightarrow \Cact(n+m-1)
\end{equation}
by the following procedure: given two cacti without spines we
re-parameterize the outside circle of the second cactus to have
length $r_i$ which is the length of the i-th circle of the first
cactus. Then glue in the second cactus by identifying the outside
circle of the second cactus with the i-th circle of the first
cactus.

Notice that this gluing differs from the one above, since now a
whole cactus and not just a lobe is re-scaled.

\subsubsection{Proposition} {\sl The gluing endows the spaces
$\Cact(n)$ with the structure of a topological operad.}

\begin{proof}
Straightforward calculation.
\end{proof}

\subsubsection{Remark}
\label{notsubop}There is an obvious map from normalized spineless
cacti to spineless cacti. This map is not a map of operads, since
the gluing procedures differ. There is however a homotopy of one
gluing to the other by moving the intersection points around the
outside circle of the cactus which is glued in, so that the two
structures of quasi-operads do agree up to homotopy. This means
that the spaces $\Cact^1$ form a homotopy associative
quasi-operad and thus the homology of this quasi-operad is an
operad. On the homology level normalized spineless cacti are thus
a sub-operad of spineless cacti and moreover, since the factors
of $\mathbb{R}^n$ are contractible this sub-operad coincides with
the homology operad of spineless cacti, as we show below.

The fact mentioned before, that the cellular chains of $\Cact^1$
form an operad \cite{del} can be seen from the discussion of the
inclusion of $\Cact^1$ into $\Cact$ mentioned above which is
explained in detail below.

\subsection{Different pictorial realizations}
\label{ckrem} As exhibited in the previous paragraph, there are
two pictorial descriptions of $\Cact^1$ and $\Cact$ given by
circles in the plane and the dual black and white planar planted
tree whose edges are marked by positive real numbers - the
lengths of the arcs. There are more pictorial realizations for
(normalized) spineless cacti, which are useful.

\subsubsection{The tree of a cactus without spines}
In the case that the configuration of circles is a cactus without
spines there is a dual tree that we can associate to it that is a
regular tree with markings that is not black and white, but is
just planar and planted.

This is done as follows. The vertices correspond to the circles.
They are labelled by the radius of the respective circle. We will
draw an edge between two vertices if the circles have a common
point and if one circle is higher than the other in the height of
the dual graph. We will label the edge by the length of the arc on
the lower circle between the intersection point and the previous
intersection point where we now also allow the length of the arcs
to be zero if these two points coincide. Here we also consider the
global zero as an intersection point.   In this procedure we give
the edges the cyclic order that is dictated by the perimeter. This
means that now the labels on the edges are in ${\Bbb R}_{\geq 0}$
with the restriction that at each vertex the label (radius) of
that vertex is strictly greater than the sum of the labels
(weights) of the incoming edges. For normalized spineless cacti
the labels on the vertices are all $1$ and can be omitted. Using
this structure we can view the space of normalized cacti as  a
sort of  ``blow up  of a configuration space''. The ``open part''
is the part with only double points. In this case, the weight on
the edges are restricted by the equations $0<\sum w_i<1$. Allowing
intersections of more than two components at a time amounts to
letting  $w_i\rightarrow 0$. In the limit $\sum w_i \rightarrow 1$
the tree is identified with the tree where the last incoming edge
is transplanted to the other vertex of the outgoing edge in such a
way that it is the next edge in the cyclic order of that vertex.
Lastly if the weight on the first edge of the root goes to zero,
the root vertex will be the other vertex of that edge.

If we do not want to use the height function of the black and white tree,
we can still define a height function via the outside circle.
Start at height zero for the root. If the perimeter hits a component for
the first time, increase the height by one and assign this height to the
component. Each time you return to a component decrease the height by one.

Given a planar planted tree whose vertices and edges are labelled in
the above fashion, it gives a prescription on how to grow a
cactus. Start at the root and draw a based loop of length given by
the label of the root. For the first edge mark the point at the
distance given by the label of the edge along the loop. Then mark
a second point by travelling the distance of the label of the
second edge and so on. Now at the next level of the tree draw a
loop based at the marked point of the previous level and again
mark points on it according to the outgoing edges. This will
produce a cactus without spines.

 Lastly, we wish to point out that  now the composition looks
 like the grafting of trees into vertices as in the
 Connes-Kreimer \cite{CK} tree operads. In fact, we have recently shown
 \cite{del} that indeed there is a cell decomposition of spineless
 normalized cacti whose cellular chains form an operad and whose
 symmetric top dimensional cells are isomorphic as an operad to
 the operad of rooted trees whose Hopf algebra is that of Connes
 and Kreimer \cite{del}.

\subsubsection{The chord diagram of a cactus} There is yet another
representation of a cactus. If one regards the outside loop, then
this can be viewed as a collection of points on an $S^1$ with an
identification of these points, plus a marked point corresponding
to the global zero. We can represent this identification scheme by
drawing one chord for each pair of points being identified as the
beginning and end of a circle this chord is oriented from the beginning
point of the lobe to the end point of the lobe. Note that one of the two segments of the
outside loop defined by the chord corresponds to the lobe.
There is a special case for the
chord diagram which is given if there is
a closed cycle of chords. This happens  if two or more
lobes intersect at the global zero. Here one can delete the first chord, if so desired,
we call this the reduced chord diagram.

The chord diagram comes equipped with a decoration of its arcs by
their length thus giving a map of $S^1_R$ to the outside circle. Here $R=\sum_i r_i$ where
the $r_i$ are the radii of the lobes. To obtain a cactus from
such a diagram, one simply has to collapse the chords.

This kind of representation is reminiscent of Kontsevich's
formalism of chord diagrams (cf.\ eg.\ \cite{BN}) as well as the
shuffle algebras and diagrams of Goncharov \cite{Go}. We wish to
point out that although the multiplication is similar to
Kontsevich's and also could be interpreted as cutting the circle
at the global zero resp.\ the local zero, it is not quite the
same. However, the exact relationship and the co-product deserve
further study.

Lastly, we  can recover the a planar rooted tree above as the dual
tree of the chord diagram. This is the dual tree on the surface
which is given by the disc whose boundary is the outside circle.
The chords on the surface then divide the disc up into chambers
--- the connected components of the complement of the chords.
The dual tree on this surface has one vertex for each  such
chamber and an edge for each pair of chambers separated by a
common chord. If the global zero lies on only one lobe the root of
the tree is the vertex of the complementary region whose boundary
includes the global zero. If there are  two components meeting at
the global zero the root of the tree is given by the vertex whose
chamber has the global as left boundary on the outside circle. In
a special case for the chord diagram which is given if there is a
closed cycle of chords, i.e.\ three or more lobes intersect at the
global zero, the root vertex will be the unique vertex inside the
closed cycle. These trees are in fact planted due to the linear
order they inherit from the embedding of the chord diagram. The
planar tree is the tree obtained from the bi-partite planted
planar tree by removing the black vertices with the exception of
the root.

\begin{figure}[ht!]
\epsfxsize = \textwidth \epsfbox{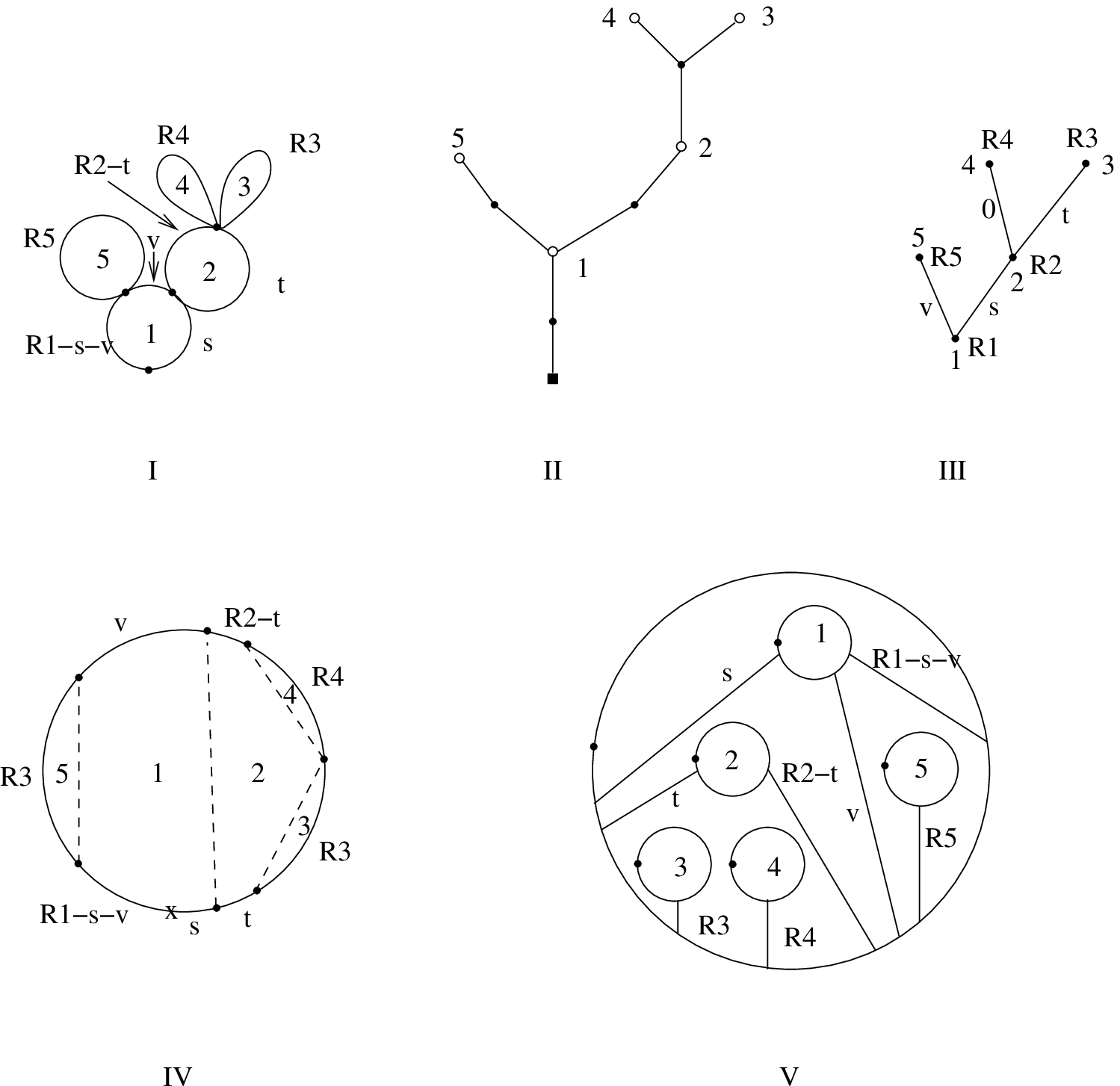}
\caption{\label{nospinesfig} I: A cactus without spines\qua II: Its
black and white tree\qua III: Its dual tree\qua IV: Its chord diagram\qua
V: Its image in $\Darc$}
\end{figure}

A representation of a cactus without spines in all possible ways
including its image in the $\Darc$ operad can be found in Figure
\ref{nospinesfig}.

\subsection{Cacti with spines}

The following definition is the original definition of cacti due
to Voronov.
\subsubsection{Definition}\cite{V}\qua  Define {\em Voronov cacti or cacti with
spines or simply cacti} in the same fashion as cacti without spines,
but without requiring that the zeros be the intersection points.

In addition, and this is key, we add a global zero/base point to
the configuration, which means that we mark a circle and a point
on that circle. The circle with the global base point will be the
root. We call the n-th component of this operad $\Cacti(n)$.

As a set $ \Cact(n) :=\{$rooted tree-like configurations of $n$
labelled $S_r^1$s in the plane$\}$/ isotopies preserving the
incidence conditions.

The perimeter or outside circle will be given by the same
procedure as \ref{perimeter} by starting at the global zero.

\subsubsection{Remark}
To define the topology we remark that the cacti with spines are as
a set in bijective correspondence to spineless cacti times a
product of $S^1$s: $\Cacti(n) \stackrel{1-1}{\longleftrightarrow}
\Cact(n)\times (S^1)^{\times n}$.

The bijection is given by mapping the underlying spineless cactus,
which is obtained by forgetting all local zeros and the induced
coordinates of the local zeros, and fixing  a coordinate on $S^1$
for every lobe which gives the length of the arc starting at the
unique intersection point with the lobe of lower height (or the
root) going counter-clockwise to the local zero.

\subsubsection{Definition}
As a topological space we set
$$
\Cacti(n):=\Cact(n)\times (S^1)^{\times n}.
$$

\subsubsection{Remark} Originally the topology was introduced by
describing that the lobes, the special points and the root can
move with the caveat that the root passes to a new component if
the intersection point of the lobe collides with the root from the
right --- just as for spineless cacti. Of course these two
descriptions are compatible. Again one can also realize
$\Cacti\subset \Darc$ and obtain the same topology as above in
this way.

\subsubsection{Gluing} We define the following operations
\begin{equation}
\circ_i: \Cacti(n) \times \Cacti(m) \rightarrow \Cacti(n+m-1)
\end{equation}
by the following procedure which differs slightly from the above:
given two cacti without spines we re-parameterize the outside
circle of the second cactus to have length $r_i$ which is the
length of the $i$-th circle of the first cactus. Then glue in the
second cactus by identifying the outside circle of the second
cactus with the $i$-th circle of the first cactus. We stress that
now the local zero of the $i$-th circle is identified with the
global zero. viz.\ the starting point of the outside circle. This
local zero need not coincide with the intersection point with the
lobe of lower height (or the global zero).

\subsubsection{Proposition}\cite{V}\qua{\sl The cacti form a topological operad.}

\subsection{Normalized cacti}
\subsubsection{Definition} We define the spaces of
{\em normalized cacti} denoted by\break $\Cacti^1(n)\subset\Cacti(n)$ to
be the subspaces of cacti with the restriction that all circles
have radius one.

As spaces
$$\Cacti^1(n) = \Cact^1(n)\times (S^1)^{\times n}.$$

\subsubsection{Glueings}
 We define the
compositions by scaling as for normalized spineless cacti and then
gluing in the second cactus into the $i$-lobe of the first, but
now using the identification of the outside circle of the second
cactus with the circle of the $i$-th lobe by matching the local
zero of the $i$-th lobe of the first cactus with the  global zero
of the second.

\subsubsection{Proposition}{\sl Together with the $\Sn$ action permuting the labels and the glueings
above normalized cacti form a topological quasi-operad.}

\begin{proof}
Straightforward computation.
\end{proof}

\subsubsection{Remark} The contents of Remark \ref{notsubop} applies
analogously in the cacti situation.

\subsubsection{Remark} There are natural forgetful morphisms from
cacti to cacti without spines forgetting all the local zeros. We
arrange the map in such a way, that the global zero becomes the
base-point of the spineless cactus. This works for the normalized
version as well. These maps are not maps of operads. The precise
relationship between the different varieties is that of a
bi-crossed product of section \ref{bicrossedsection}, see Theorem \ref{bicrossed} below.

There is, however, an embedding of spineless cacti into cacti as a
suboperad by considering the global zero to be the zero of the
root and by making the first intersection point at which the
perimeter reaches a lobe of the cactus the local zero of that
circle.

\subsection{Different pictorial realizations}
\subsubsection{The tree of a cactus}

The missing information of a cactus without spines relative to
a cactus proper is the location of the local zeros. We just add this
information as a second label on each vertex.
Notice that the local zero of the root component then need not be the global
zero. The label we associate to the root is the position of
the local zero with respect to the global zero.

\subsubsection{The chord diagram} The chord diagram of a cactus
again is the chord diagram of a cactus without spines, where the
location of the spines is additionally marked on the $S^1$. There
is a choice if the local zero coincides with an intersection
point. Just to fix notation we will mark the first occurrence of
the endpoint of a chord, where first means in the natural
orientation starting at the global zero.

\begin{figure}[ht!]
\epsfxsize = \textwidth \epsfbox{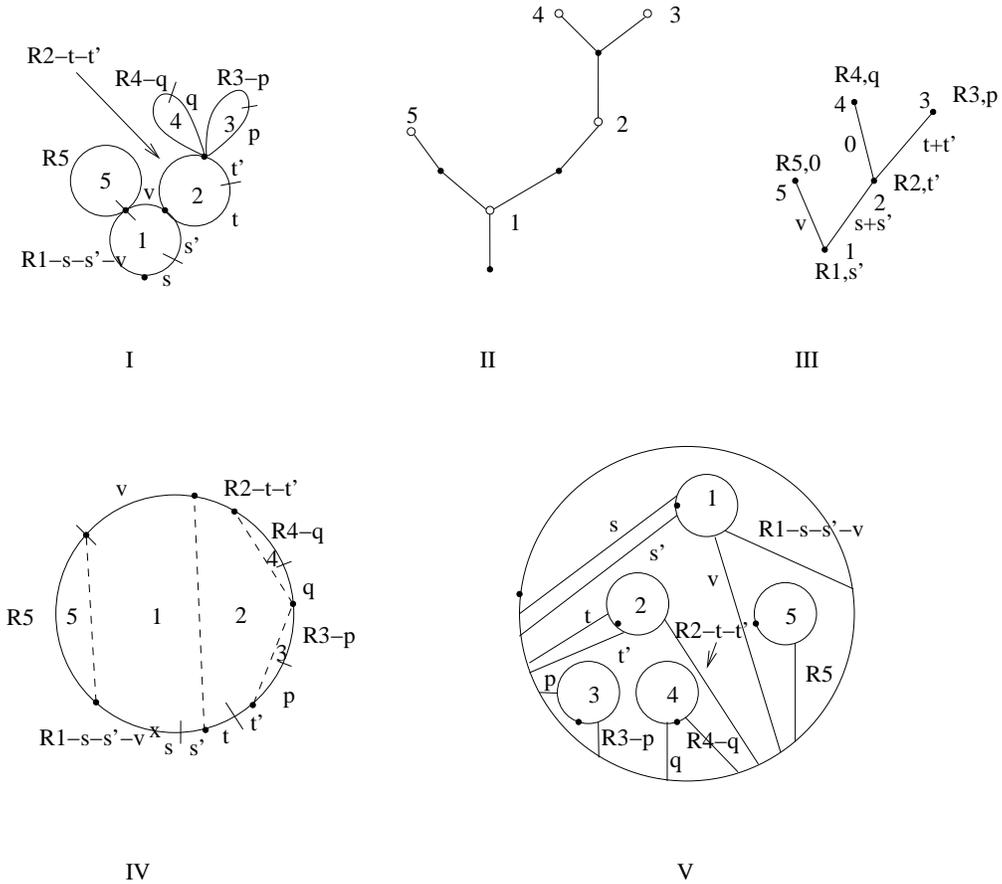}
\caption{\label{spinesfig} I: A cactus (with spines)\qua II: Its black
and white tree\qua III: Its dual tree\qua IV: Its chord diagram\qua V: Its
image in $\Darc$}
\end{figure}

A representation of a cactus (with spines) in all possible ways
including its image in $\Darc$ can be found in Figure
\ref{spinesfig}.

\section{Spineless cacti and the little discs operad}

In this section, we will show that $\Cact$ is an $E_2$ operad using the recognition
principle of Fiedorowicz \cite{fiedorrecog}. To assure the needed assumptions are
met we mimic  the construction of \cite{fiedorbar} which shows that the universal
covers of the little discs operad naturally form a $B_{\infty}$ operad.

\subsection{The $E_1$ structure}

\subsubsection{Definition} A spineless corolla cactus (SCC) is a
spineless cactus whose points of intersection all coincide with
the global zero.

Since the condition of being an SCC is preserved when composing two spineless corolla cacti:
\subsubsection{Lemma}
{\sl Spineless Corolla Cacti are a suboperad of spineless Cacti.}

We define $\Corol(n)\subset\Cact(n)$ to be the subset of
spineless corolla cacti and
denote the operad constituted by the $\Corol(n)$ with the
permutation action of $\Sn$ and the induced gluing by
$\Corol(n)$.

\subsubsection{Theorem}{\sl The suboperad $\Corol(n)$ of corolla
cacti is an $E_1$ operad.}

\begin{proof}
We will use the recognition principle of Boardman-Vogt \cite{BV}.
First notice that we have a free action of $\Sn$. If the lobes are
all grafted together at the root then the only parameters are the
sizes of the lobes. These sizes together with the labelling fixes
a unique spineless corolla. Two spineless corollas lie in the same
path component if and only if the sequence of the labels of the
lobes as read off from the outside circle agree. Thus
$\Corol(n)=\amalg_{\s\in\Sn} \mathbb{R}^n_{>0}$. And thus each
path component is contractible and thus the action of $\Sn$ is
free and transitive on $\pi_0(\Corol(n))$.
\end{proof}

\subsubsection{Remark}
The Theorem above has immediate applications to operads built from
moduli spaces (see Appendix B) giving them an $A_{\infty}$-structure.

\subsubsection{Corollary}
\label{eone}{\sl The operad of the decorated moduli space of bordered,
punctured surfaces with marked points on the boundary $\widetilde
M^s_{g,r}$ which is proper homotopy equivalent to
 $\Arc_{\#}$ contains an $E_1$ operad. Thus so does the operad
 $\Arc$.

 The same is true for the
 operad of the moduli spaces $M_{g,n}^n$ of genus $g$ curves with $n$
 punctures and a choice of tangent vector at each puncture
 and its restriction to genus
 $0$.

 Finally the spaces $M_{g,n}$ form a partial operad which is an
 $E_1$ operad.}

\begin{proof}
By the Appendix B there is an operad map $\Cacti
\rightarrow \Arc_{\#}\subset \Arc$ which is an equivalence
onto its image. Furthermore it is shown that $\Arc_{\#}$ is proper
homotopy equivalent to the mentioned moduli space in \cite{P2}.
This establishes the first part.

The second claim follows from the identification of the suboperad
of bordered surfaces with marked points on the boundary and no
further punctures $\Arc^0_{\#}$ with $M_{g,n}^n$ via marked ribbon
graphs \cite{hoch}.

The last statement comes from the fact that the SCCs are ribbon
graphs and as such index cells of $M_{g,n}$. The operad structure
of SCCs thus defines a partial operad structure on $M_{g,n}$.
\end{proof}

This also means that on the chain level algebras over these operads
will be $A_{\infty}$ operads.

\subsection{The $E_2$ structure}

The main result of this section is the following.

\subsubsection{Theorem}
\label{cactdisc} {\sl $\Cact$ is an $E_2$ operad.}

We will use the recognition principle of Fiedorowicz
\cite{fiedorrecog} to prove this theorem (see also \cite{SW}). For
this one needs the notion of a braid operad, which is 
given by replacing the symmetric groups in the definition of
operads by braid groups (see \cite{fiedorbar} or \cite{SW}).

\subsubsection{Definition}\cite{fiedorbar}\qua
A collection $B(n)$ a $B_{\infty}$ operad if the $B(n)$ form a
braid operad in the sense of \cite{fiedorbar} with the properties

\begin{itemize}
\item[i)] the spaces $B(n)$ are contractible.

\item[ii)] The braid group action on each $B(n)$ is free.
\end{itemize}

\subsubsection{Proposition}\cite{fiedorrecog}\qua
{\sl An operad $\mathcal{A}$ is an $E_2$ operad if and only if each
space $\mathcal{A}(k)$ is connected and the collection of covering
spaces $\{\tilde A(k)\}$ form a $B_{\infty}$ operad.} \endproof

Adapting the proof of \cite{fiedorbar} that the universal covers
of the little disc operad form a $B_{\infty}$ operad one arrives
at the following proposition, which is essentially contained in
\cite{fiedorbar} and in spirit in \cite{MS}.

Let $\t_i \in \Sn$ denote the transposition which transposes $i$ and $i+1$.

\subsubsection{Proposition} \label{MSprop} {\sl Suppose we are given an
operad $D(n)$ with the properties:

\begin{itemize}
\item[\rm i)] The $\Sn$ action on each $D(n)$ is free.

\item[\rm ii)] $D$ affords a morphism of non-$\Sigma$ operads $I$ :
$C_1 \rightarrow D$ where $C_1$ is an $E_1$ operad.

\item[\rm iii)] $D(n)/\Sn$ is a $K(Br_n,1)$ for the braid group $Br_n$
where the braid action covers the symmetric group action.

\item[\rm iv)] The spaces are $D(n)$ are homotopy equivalent to
CW complexes.
\end{itemize}
 Then the collection of universal covers $\tilde D(n)$
 is a $B_{\infty}$ operad and hence $D$  is equivalent as an
operad to $C_2$, the little 2-cubes operad.}

\begin{proof}
Let $p:\tilde D(n) \rightarrow D(n)$ be the universal cover. We
have to show that the spaces $\tilde D(n)$ form a braid operad and
that they are contractible. The latter fact is true since by iii)
the spaces $\tilde D(n)$ are weakly contractible and by assumption
iv) the $D(n)$ are homotopic to a CW complex, so that the $\tilde
D(n)$ are indeed contractible. For each $n$ choose a component of
$p^{-1}(I(C_1(n))$ which we call $\tilde C_1$.

The $\tilde C_1$ allow to lift the operad composition maps by
letting $\tilde \gamma$
$$
\begin{CD}
\tilde D(k)\times \tilde D(j_1) \dots \times \tilde D(j_k)@>\tilde
\gamma>> \tilde D(j_1+\dots +j_k)\\
@VVpV @VVpV\\
 D(k)\times  D(j_1) \dots \times D(j_k)@>
\gamma>> D(j_1+\dots +j_k)\\
\end{CD}
$$
be the unique lift which takes $\tilde C_1(k)\times \tilde
C_1(j_1) \dots \times \tilde C_1(j_k)$ to $\tilde C_1(j_1+\dots
+j_k)$. To write out the braid action fix a point $c_n \in C_1(n)$
and for each $i$ a path $\a_i$ from $I(c_n)$ to $\t_iI(c_n)$ which
lifts a non-null homotopic path of $D(n)/\Sn$. Notice that these
satisfy the conditions that for each $n$ and $i$ the paths $
\t_i\t_{i+1}(\a_i) \cdot \t_i(\a_{i+1}) \cdot \a_i$ and
$\t_{i+1}\t_i(\a_{i+1}) \cdot \t_{i+1}(\a_i) \cdot \a_{i+1}$ are
path homotopic (where $\cdot$ denotes concatenation of paths) due
to condition iii).
 The explicit paths $\t_i$  then provide the $Br_n$ action on $\tilde
D$(n) again by using $\tilde C_1(n)$ as ``base-points'' to lift the $\Sn$ action.
 It is now a
straightforward computation that the compositions $\tilde \gamma$
and the braid group action define a braid operad. Furthermore the
braid group actions are free by iii) and thus the $\tilde D(n)$
form a $B_{\infty}$ operad.
\end{proof}

\begin{proof}[Proof of Theorem \ref{cactdisc}]
As announced we will check the conditions of Prop\-osition
\ref{MSprop}. In our case the operad $D$ will be the operad of
spineless cacti $\Cact$. The condition i) is obvious, since $\Sn$
acts freely on the labels. We showed above that $\Corol$ is an
$E_1$ operad which is a suboperad of $\Cact$. This establishes
ii). The condition iii) follows from Proposition \ref{piprop}
below. Lastly, the condition iv) follows from the definition of the spaces
$\Cact(n)=\Cact^1(n) \times \mathbb{R}_{>0}^n$.
\end{proof}

\subsection{The forgetful quasi-fibration}
This section is devoted to showing that the spaces $\Cact(n)/\Sn$ are
$K(Br_n,1)$.

\subsubsection{Definition} The completed chord diagram of a cactus $c$ without spines
is the topological space obtained as follows. Cut the outside
circle at the global zero, mark the two endpoints and add a chord
$a_z$ between them. If a chord started at the global zero, then
the new starting point will be the right endpoint of  $a_z$, if
ended on the root then the new endpoint will be the left endpoint
of $a_z$.

Identify each marked point  (that is the added endpoints of $a_z$ and the
endpoints of the chords) of the circle with a 0-simplex and each
arc connecting two marked points with a  1-simplex joining the two
0-simplices. Now for any sequence of chords connecting $k$
points of the outside circle glue in a $k-1$ simplex, by
identifying the vertices of the simplex with these points. For any
chord including $a_z$ this means that the sub-one-simplex given by
the two endpoints of the chord can be identified with the chord.
We let the diagram have the co-induced topology.

\begin{figure}[ht!]
\cl{\epsfxsize = .8\textwidth
\epsfbox{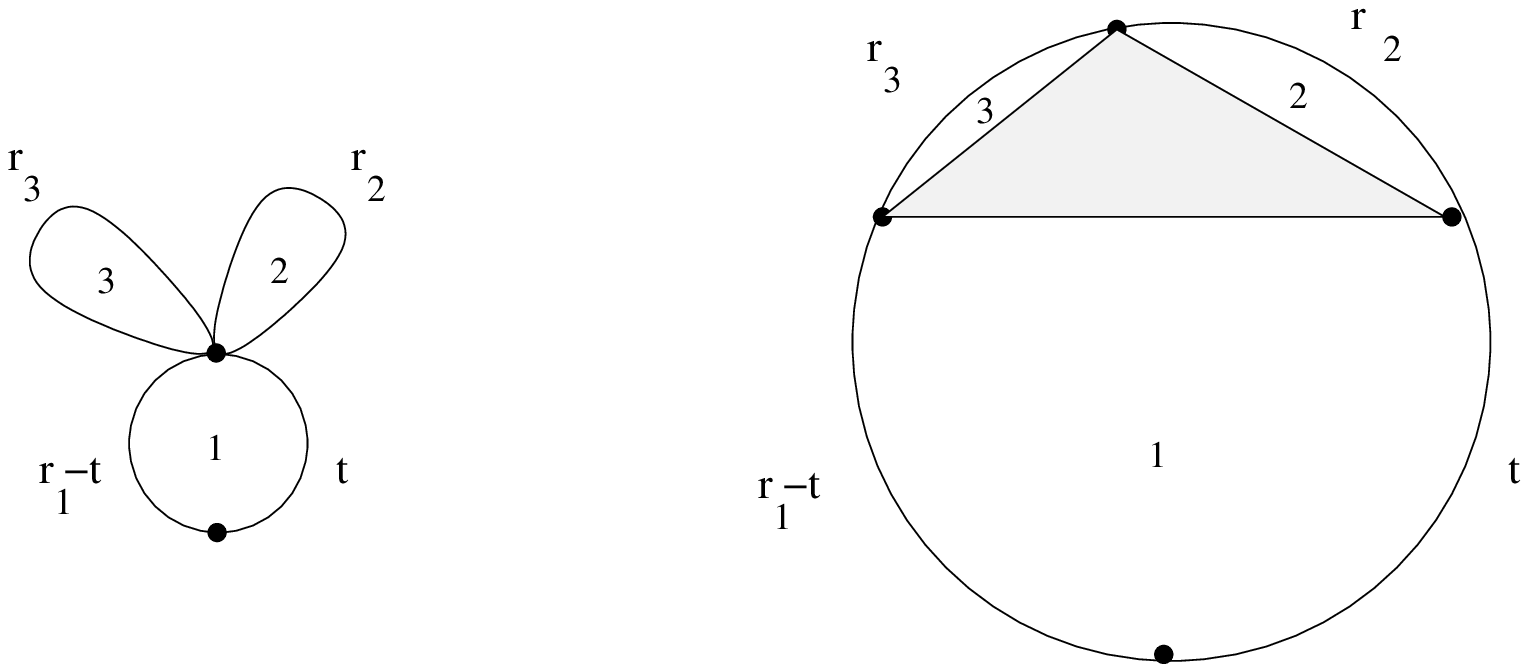}}
\caption{\label{compchord} A cactus without spines and its completed chord diagram}
\end{figure}

For an example of a completed diagram, see Figure \ref{compchord}.

\subsubsection{Definition} We define the spine of a completed chord
diagram to be the following subspace. For each maximal k-simplex, fix the
barycenter. First connect the barycenter to all the vertices of
the simplex by a straight line, then connect the vertices of the
simplices by the arcs of the outside circle to obtain the spine.

\subsubsection{Lemma}
\label{retract}
{\sl A completed chord diagram of a cactus $c$ is
homotopy equivalent to its spine which is homotopy equivalent to
the image of $c$.}

\begin{proof}
By retracting to the spine, the first claim follows. For the
second claim in one direction we contract of the straight lines of
the spine to retrieve the underlying cactus. For the homotopy
inverse identify the vertices of the cactus with the barycenters.
Each arc of the cactus $a$ corresponds to a unique arc $a'$ on the
outside circle of the chord diagram. Let $v_1,v_2$ be the
starting- and the endpoint of the directed arc. Now map each arc
$a$ of the cactus between two vertices to the path between the
barycenters representing these vertices which first goes from the
barycenter to the vertex $v_1$, then along the arc $a'$, and
finally from $v_2$ to the second barycenter. .
\end{proof}

Now we will consider the surjective map $p_{n+1}:\Cact(n+1)
\rightarrow \Cact(n)$, which contracts the $n+1$-st lobe of the
cactus. We call the image of the contracted lobe the marked point.
If the root happens to lie on the component $n+1$ then the root
after the contraction is fixed to be the marked point.

\subsubsection{Forgetful maps}
Define a map
$$p^T:\wlbptree(n+1)\rightarrow \wlbptree(n)$$
 by mapping a labelled
tree $\t\in \wlbptree(n+1)$ to  $p^T(\t)\in \wlbptree(n)$ which is
the tree obtained from $\t$ by  coloring the vertex $v_{n+1}$
labelled by $n+1$ black, forgetting the label and contracting all
the edges incident to this vertex. If the image of the vertex
$v_{n+1}$ under the contraction only has one adjacent edge, we
also delete this vertex and this edge to define $p^T(\t)$.

This induces projection maps $p^{\Delta(\t)}:\Delta(\t)\rightarrow
\Delta(p^T(\t))$ projecting to the product of the first $n$
simplices, i.e.\  forgetting the coordinates of the flags of $v$.
Formally, let $E(\t)$ be the edges of $\t$, and $E(v)$ be the
edges incident to $v$. We map the point with coordinates $(x_e),
e\in E(\t)$ in $\Delta(\t)$ to the point with coordinates
$(x'_e=x_e), e\in E(p^T(\t))$ in $\Delta(p^T(\t))$, where we
identified the non-contracted edges of $\t$ with those of
$p^T(\t)$.

Now we define a map $p':\Cact^1(n+1)\rightarrow \Cact^1(n)$ as
follows. For $c'\in \Cact^1$ let $\t$ be its topological type. Set
$$p'(c'):=e_{p^T(\t)}\circ p^{\Delta(\t)} \circ \dot e_{\t}^{-1}(c)$$
 Finally let
$c=(c,(r_1,\dots,r_{n+1})) \in \Cact_{n+1}$. We define
$$p(c',(r_1,\dots,r_{n+1}))=(p'(c'),(r_1,\dots, r_n))$$
This defines the map $p:\Cact(n+1)\rightarrow \Cact(n)$ mentioned
above.

\subsubsection{Proposition}
\label{fiberchord}{\sl The fiber of the map $p$ over a spineless
cactus $c$ is homeomorphic to the completed chord diagram of $c$
times $\mathbb{R}_{>0}$. The fiber of the map $p'$ over a
normalized spineless cactus $c$ is homeomorphic to the completed
chord diagram of $c$.}

\begin{proof}
This follows directly from the above description of the map $p'$
as projecting out the simplex for the vertex $v_{n+1}$ labelled by
$n+1$.

 A detailed description is as follows: Fix a cactus $c\in
\Cact(n)$. The fiber over it can then be characterized in the
following way. First there is a factor of $\mathbb{R}_{>0}$ which
fixes the radius of the cactus. Then there is an interval which
parameterizes the cacti where the lobe $n+1$ has no lobe above it
and does not contain the root. The parametrization is via a marked
point on the outside circle. In the case that the root is only on
one lobe this interval is glued together with another interval to
form a circle. The second interval parameterizes the pre-images
obtained by gluing the $n+1$st lobe to the root and then moving
the root around that lobe. If the root is moved all the way
around, the limit is the same configuration as the one in which
the lobe has been moved all the way around the outside circle.
Also keeping the root at the intersection point is the same
configuration as the initial point of the first interval. Thus we
obtain a circle glued from two intervals.

If the pre-image is such that the $n+1$st lobe has higher lobes
attached to it, then the marked point is  necessarily a point of
intersection on the cactus. First assume that this intersection
point is not the global zero. This means that we can blow up this
intersection point to a circle and arrange the lobes attached to
it keeping their order according to the outside circle. Such a
configuration is determined by the length of the arcs between the
attached circles. These lengths add up to the total radius and
thus are parameterized by points in a $k$-simplex, if the number
of lobes meeting at this point is $k+1$. At the vertices of this
simplex each lobe is again attached to the common intersection
point. These configurations coincide with the points of the chord
diagram to which the chords are attached. In the case that the
global zero is at the intersection of $k>1$ lobes, then we again
``blow up'' the global zero to a $k$-simplex. One of the edges of
this simplex is the interval in which the $n+1$ lobe is attached
to the global zero and the global zero is moved around this lobe
as discussed above and identified with the arc replacing the
global zero when completing the chord diagram. The $k-2$ vertices
which are not on this edge are then identified with the $k-2$
vertices of the closed sequence of chords excluding the global
zero.
\end{proof}

\subsubsection{Corollary}
{\sl The fiber of the map $p$ over a spineless cactus $c$ is homotopy
equivalent to the image of the cactus $c$ in $\mathbb{R}^2$ and is
thus homotopy equivalent to a bouquet of $n$ circles $\bigvee_n
S^1$.}

\begin{proof}
The first equivalence follows from Lemma \ref{retract}. The second
equivalence is straightforward.
\end{proof}

\subsubsection{Remark} 
Let $\UCact(n)$ be the set obtained from $\Cact(n+1)$ by contracting the simplices of
the completed chord diagrams in each fiber of the map
$p:\Cact(n+1)\rightarrow\Cact(n)$. Let $\rho_{n+1}:\Cact(n+1)\rightarrow \UCact(n)$ 
be the induced surjection and endow $\UCact(n)$ with the quotient topology.
The map $p$ factors though $\rho$, that is $p_{n+1}=\tilde p_{n} \circ \rho_{n+1}$ where 
$\tilde p_{n}:\UCact(n)\rightarrow \Cact(n)$ is the universal map whose fiber over
a spineless cactus $c$ is the image of that spineless cactus and whose total space is
$\UCact(n)=\bigcup_{c\in \Cact(n)}Im(c)$. A point in this space is a cactus together with an additional marked point on the cactus. The map $\tilde p$ forgets this point.

Then $p:\Cact(n+1)\rightarrow\Cact(n)$ 
is fiberwise homotopy equivalent to the universal map $\tilde p_n:UCact(n)\rightarrow \Cact(n)$.

\subsubsection{Remark} The maps $p$ and $p'$ are not fibrations.
We will show that they are,  however,  quasi-fibrations.

\subsubsection{Definition}
Let  $c,c' \in \Cact^1$ and $\t,\t'$ be their topological types.
We say that $c'$ can be derived from $c$ and also that $\t'$ can be
derived from $\t$ if $e_{\t'}(\Delta(\t'))\subset
e_{\t}(\Delta(\t))\subset \Cact^1$. Here $\D(\t)$ is the product of simplices as defined in
equation (\ref{deltataudef}).

If the inclusion is proper, we say that $c'$ is a degeneration of
$c$ and also say $\t'$ is a degeneration of $\t$.

\subsubsection{Remark}
The notion of degeneration induces a partial order on $\wlbptree$
where $\t'\prec\t$ if $\t'$ is a degeneration of $\t$.

\subsubsection{Definition}
If there is a $\t''$ s.t.\ $c,c' \in e_{\t''}(\Delta(\t''))$ we
say $c,c'$ share the common type $\t''$ and also say that $\t,\t'$
share the common type $\t''$.

In  case $c$ and $c'$ share a common type $\t''$, we let
$d_{\t''}(c,c')$ be the distance between their lifts into
$\Delta(\t'')$.

\subsubsection{Definition}
For $c\in \Cact(n)$, $\t$ with $c\in e_{\t}(\Delta(\t))$ and
$\eps>0$ we define
\begin{eqnarray}
U(c,\eps,\t)&:=&\{c'\in  e_{\t}(\Delta(\t))|d_{\t}(c,c')<\eps\}
\mbox{ and }\nn\\
U(c,\eps)&:=&\bigcup_{\t: c\in e_{\t}(\Delta(\t))} U(c,\eps,\t)
\end{eqnarray}
It is clear that the $U(c,\eps)$ are open.

We call $\eps$ small for $c$ if $c'\in U(c,\eps)$ implies that $c$
is a degeneration of  $c'$.

\subsubsection{Remark}
The set of small $\eps$ for a fixed $c$ is non-empty. For instance
if $A$ is the set of arcs of $c$, any
$\eps<\frac{1}{2}\min(|a|:|a|\neq0)$ will do, since one cannot
move the root or a lobe more than the length of any arc and hence
cannot create new degenerations without going beyond the distance
$\eps$.

\subsubsection{Lemma}$\phantom{99}$
\label{ausgezeichnet}
{\sl\begin{itemize}
\item[\rm i)] The sets $U(c,\eps)$ with $\eps$ small for $c$ are open
and contractible.

\item[\rm ii)] The sets $U(c,\eps)$ with $c\in \Cact^1(n)$ cover
$\Cact^1(n)$.

\item[\rm iii)] If $c''\in U(c,\eps)\cap U(c',\eps')$ then there
exists an $\eps''$ s.t.\ $c''\in U(c'',\eps'')\subset
U(c,\eps)\cap U(c',\eps')$.
\end{itemize}}

\begin{proof} The fact that these sets are open and cover is immediate.
For the contraction we define $h:U(c,\eps) \times I\rightarrow
U(c,\eps)$ as follows: for $c' \in U(c,\eps) $  with topological
type $\t'$, we set $c'(t)$ to be the image of the point of
$\Delta(\t')$ which is at distance $\eps-t/\eps$ from the point
corresponding to $c$ in $\Delta(\t')$ along the unique line
joining these two points. This map is easily seen to be continuous
and contracts $U(c,\eps)$ onto $c$.

For part ii) let $\t,\t',\t''$ be the respective topological
types. If  $d_{\t''}(c,c'')=d_1$ and $d_{\t''}(c',c'')=d_2$ fix
any  $\eps''<\min(\sqrt{\eps^2-d_1^2},\sqrt{\eps^{\prime
2}-d_2^2})$. Then the inclusion follows from the fact that $\prec$
is a partial order, i.e.\ if $\t''\preceq \t'''$ then also
$\t\preceq \t'''$ and $\t'\preceq \t'''$ since $\t\preceq\t''$ and
$\t'\preceq \t''$.
\end{proof}

\subsubsection{Remark} The contraction
simultaneously contracts all those arcs of $c'$ which do not
appear in $c$, i.e.\ those which correspond to the edges of $\t'$
which are contracted to obtain $\t$.

\subsubsection{Lemma}
\label{localcheck}{\sl The pair $(p'^{-1}(U(c,\eps)), p'^{-1}(c))$ is
homotopy equivalent to $(p'^{-1}(c), p'^{-1}(c))$.}

\begin{proof}
We define the homotopy
$$H:(p'^{-1}(U(c,\eps)), p'^{-1}(c))\times I
\rightarrow (p'^{-1}(U(c,\eps)),p'^{-1}(c))$$
as follows. Given
$\hat c'\in p'^{-1}(U(c,\eps))$ write it as the tuple $(c',
ch(c'))$ where $c'=p'(\hat c')$ and $ch(\hat c')$ is the point in
the fiber over $c'$. By Proposition \ref{fiberchord} $ch(\hat c')$
is a unique point of the completed chord diagram $Chord(c')$ of
$c'$. We let $H(t)(c')=( c'(t),ch(c')(t))$ where $ c'(t):=h(
c',t)$ is the cactus as in Lemma \ref{ausgezeichnet} and $ch(\hat
c')(t)\in Chord(c'(t))$ is defined as follows. First notice that
during the homotopy $h$ the topological type $\t'$ of $c'$ does
not change as long as $t\neq 1$ and therefore there are natural
homeomorphisms $h_{chord}(t):Chord( c')\rightarrow Chord( c'(t))$
obtained by a homogeneous re-scaling of the arcs of the outside
circle by factors $x_a(c)(t)/x_a(c)$ -where again the $x_a$ are
the coordinates in $\Delta(\t)$. For $t\neq 1$ we set $ch(\hat
c')(t):=h_{chord}(t)(ch(\hat c'))$.

In order to extend to $t=1$ notice that the chords of $Chord(c')$
and those of $Chord(c)$ are in 1-1 correspondence as 1-simplices,
as they correspond to the labels 1 through n. Therefore the
simplices of $Chord(c')$ uniquely correspond to faces of the
simplices  or simplices of $Chord(c)$. We let
$h_{chord}(1):Chord(c')\rightarrow Chord(c)$ be the map that first
contracts the arcs of the outside circle which are indexed by arcs
$a$ with $x_a(c)=0$ and then identifies the result of this
contraction with a subset of $Chord(c)$ by identifying the arcs of
the outside circle with the same labels and identifying the
simplices of $Chord(c')$ with the respective faces of $Chord(c)$.
Finally set $ch(\hat c')(1)=h_{chord}(1)(ch(\hat c'))$. It is now
easy to check that the defined map is indeed a homotopy.
\end{proof}

\subsubsection{Remark}
The effect of the contraction above  is to contract the arcs not
belonging to $c$ while keeping the lobe $n+1$ in its relative place.

\subsubsection{Proposition}
\label{quasifibprop}
{\sl $p':\Cact^1(n+1) \rightarrow \Cact^1(n)$ and
$p:\Cact(n+1) \rightarrow \Cact(n)$ are quasi-fibrations.}

\begin{proof}
First let's handle $p'$: $p'|p^{\prime -1}(U(c,\eps))$ is a
quasi-fibration by the Lemma \ref{localcheck} above. This fact
together with Lemma \ref{ausgezeichnet} shows that the conditions
of the Dold-Thom criterium \cite{DT}[Satz 2.2] are met and hence
$p'$ is a quasi-fibration. Now fix some base-point
$c=(c',\vec{r})$, then the claim follows from the following
equalities:
\begin{multline}
\pi_i(\Cact(n+1),p^{-1}(c))=\pi_i(\Cact^1(n+1)\times
\mathbb{R}^{n+1}_{>0},p^{-1}(c',\vec{r}))\\
=\pi_i(\Cact^1(n+1),p'^{-1}(c'))=\pi_i(\Cact^1(n),c')=\pi_i(\Cact(n),(c,\vec{r}))
\end{multline} where the first equality holds by definition, the
second holds since the pair $(\Cact^1(n+1)\times
\mathbb{R}^{n+1}_{>0},p^{-1}(c',\vec{r}))$ is homotopy equivalent
to the pair $(\Cact^1(n+1),p'^{-1}(c'))$ by contracting the
factors $\mathbb{R}_{>0}$ to the point $1$, the third equation
holds, since $p'$ is a quasi-fibration and finally the last
equation holds since again by contraction of the factors
$\mathbb{R}_{>0}$  the pair $(\Cact^1(n),c')$ is homotopy
equivalent to the pair $(\Cact(n),(c,\vec{r}))$.
\end{proof}

\subsubsection{Proposition}
\label{piprop} {\sl The spaces $\Cact(n)$ are $K(PBr_n,1)$ spaces and
the spaces $\Cact(n)/\Sn$ are $K(Br,1)$ spaces.}

\begin{proof}

Since by Proposition \ref{quasifibprop} $p$ is a quasi-fibration,
we have the long exact sequence of homotopy groups \cite{DT}
$$
\rightarrow \pi_{i+1}(\Cact(n))\rightarrow \pi_i(\bigvee_n S^1)
\rightarrow \pi_i(\Cact(n+1)) \rightarrow \pi_i(\Cact(n))
$$
where we inserted $\pi_i(p^{-1}(c))=\pi_i(\bigvee_n S^1)$.

The fibration $p$ admits a section, for instance attaching the
(n+1)st lobe at the root and letting the root lie on the new
(n+1)st lobe. Thus the long exact sequence for this
quasi-fibration splits.

First notice that  since $\Cact(1)=*, \Cact(2)=S^1$ by induction
$\pi_i(\Cact(n))=0$ for $k\geq 2$ and so also $\pi_i(\Cact(n)/\Sn)=0$ for $k\geq 2$  .

For the first homotopy group, we fix some data. Choose the
spineless corolla cactus with radii all equal to one and labelling
$1,2, \dots,n$ as the base point $c_n$ of $\Cact(n)$ and choose
the paths $\a_i$ as indicated in Figure \ref{alphai} which makes
the braid action explicit.

\begin{figure}[ht!]
\epsfxsize = \textwidth \epsfbox{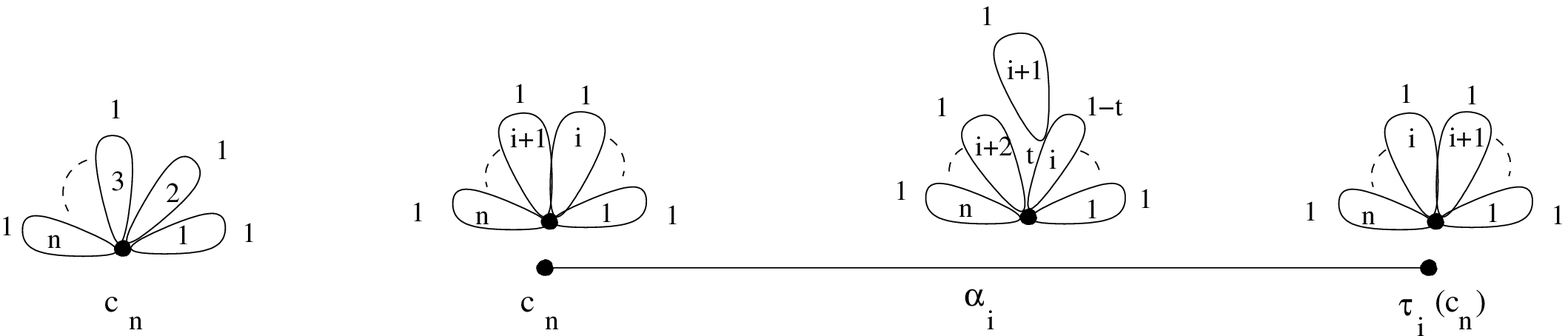}
\caption{\label{alphai}The point $c_n$ and the path $\a_i$}
\end{figure}

Note that the braid condition that for each $n$ and $i$ the paths
$ \t_i\t_{i+1}(\a_i) \cdot \t_i(\a_{i+1}) \cdot \a_i$ and
$\t_{i+1}\t_i(\a_{i+1}) \cdot \t_{i+1}(\a_i) \cdot \a_{i+1}$ are
path homotopic (where $\cdot$ denotes concatenation of paths) is
verified explicitly in Figure \ref{cactbraid} were we have only
drawn the relevant three lobes and indicated the other lobes by
dots.
\begin{figure}[ht!]
\epsfxsize = \textwidth \epsfbox{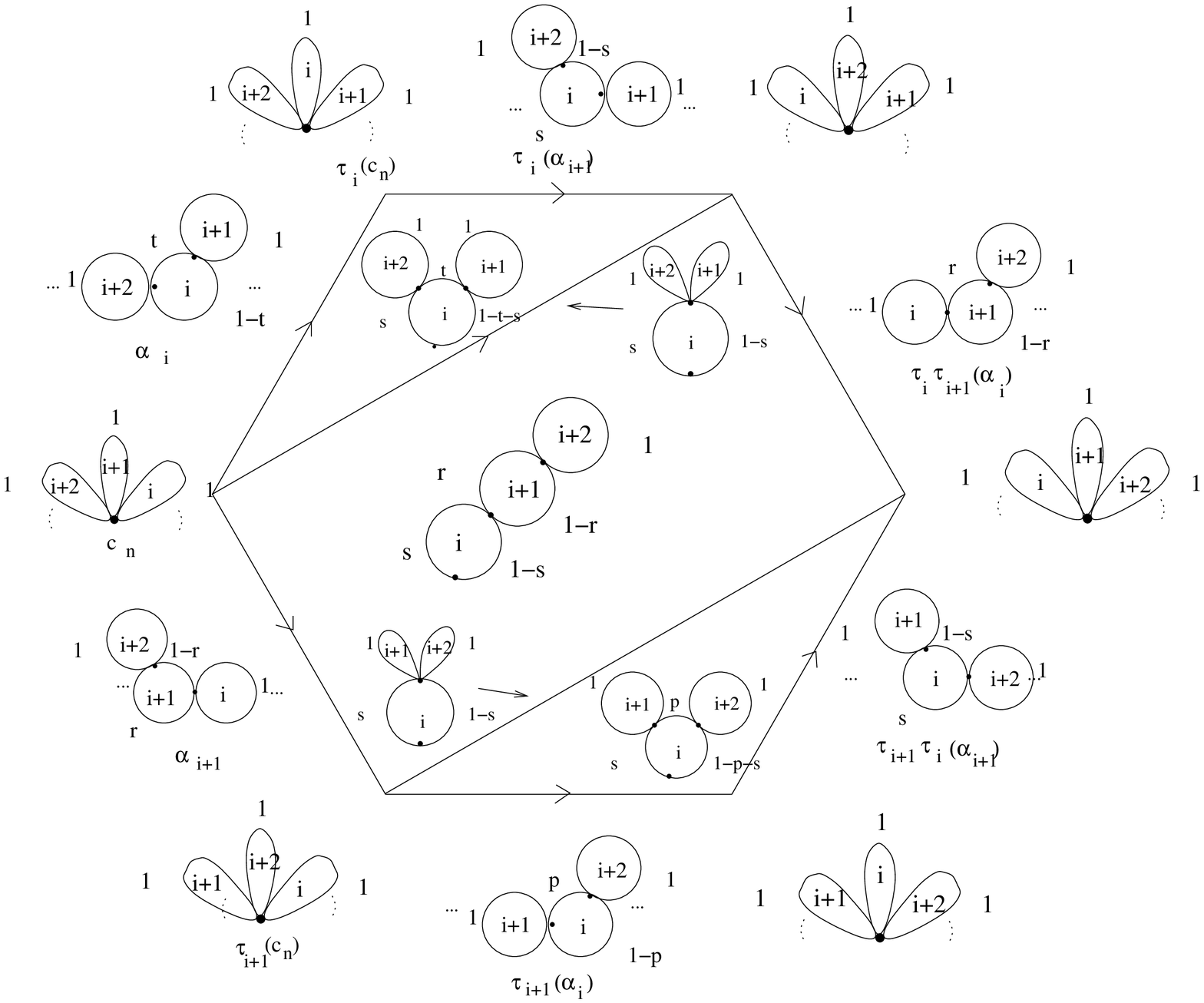}
\caption{\label{cactbraid}The braid homotopy}
\end{figure}

Now we proceed by induction on $n$ assuming
$\pi_1(\Cact(n))=PBr_n$ where $PBr_n$ is the pure braid group and
$\pi_1(\Cact(n)/\Sn)\simeq Br_n$ with the explicit maps
$Br_n\rightarrow \pi_1(\Cact(n)/\Sn)$ is given by $b_i\mapsto
\a_i$. Notice $\pi_1(\Cact(2))\simeq\mathbb{Z}=PBr_2$ generated by
$\t_1(\a_1)\cdot\a_1$ and $\pi_1(\Cact(2)/\mathbb{S}_2)\simeq
Br_2\simeq \mathbb{Z}$ generated by $\a_1$.

First we treat $\pi_1(\Cact(n+1))$. For this consider the
following diagram:
$$
\begin{CD}
1 @>>> \pi_1(\bigvee_n S^1) @>>> PBr_{n+1} @>>> PBr_n @>>>1\\
@. @| @VVV  @VVV\\
1 @>>> \pi_1(\bigvee_n S^1) @>>>
\pi_1(\Cact(n+1)) @>>> \pi_1(\Cact(n)) @>>>1\\
\end{CD}
$$
Here the second line follows from the long exact sequence and the
first line is a classic fact. It follows e.g.\ from regarding the
long exact sequence for the forgetful map between the
configuration spaces of $n+1$ and  $n$ ordered points in
$\mathbb{R}^2$ which forgets the $n+1$st point. By induction we
know the map $Br_n \rightarrow \pi_1(\Cact_n/\Sn)$ sending the
generator $b_i$ which maps to the transposition $\t_i$ in $\Sn$ to
$[\a_i]\in \Cact(n)/\Sn$ is an isomorphism. The right down arrow
is its restriction to $PBr_n$ which is the isomorphism sending the
generators $\xi_{ij}=b_i b_{i+1}\dots
b_{j-1}b_j^2b_{j-1}^{-1}\dots b_i^{-1}:1\leq i<j \leq n$  of
$PBr_n$ to the class $[\a_{ij}]\in \pi_1(\Cact(n))$ where
$\a_{ij}$ is the closed path $\t_{i+1}\dots
\t_{j-1}\t_j^2\t_{j-1}^{-1}\dots \t_i^{-1}(\a_i)\cdot \ldots
\cdot\t_j\dots\t_{i}(\a_{j})\cdot
\t_{j-1}\dots\t_{i}(\a_{j})\cdot\ldots
\cdot\t_i(\a_{i+1})\a_i^{-1}$. Choosing the base point of the
fiber to be the cactus $c_n$ we see that the generators of
$\pi_1(Chord(c_n))$ can be identified with the paths $\a_{in+1}$
hence identifying the left isomorphism as $\xi_{in+1}$ maps to
$[\a_{in+1}]$. Hence we have a diagram of group extensions and the
middle arrow which sends $\xi_{i,j}$ to $[\a_{ij}]$ for $1\leq
i<j\leq n+1$ is also an isomorphism.

For the fact that $\pi_1(Cact(n+1)/\Snn)=Br_{n+1}$ consider the
following diagram of group extensions
$$
\begin{CD}
1 @>>> PBr_{n+1} @>>> Br_{n+1} @>>> \Snn @>>>1\\
@. @VVV @VVV  @|\\
1 @>>>  \pi_1(\Cact(n+1)) @>>> \pi_1(\Cact(n+1)/\Snn)@>>>\Snn @>>>1\\
\end{CD}
$$
where the left arrow was shown to be an isomorphism above and the
commutativity follows from the explicit mapping of the upper to
the lower row which sends $b_i$ to the respective class of $\a_i$.
Hence the middle arrow is an isomorphism, proving the claim.
\end{proof}

\subsubsection{Remark}
In our proof, we chose the $E_1$ operad of spineless corollas. If
we want to use little intervals there is an obvious map by
assigning the lengths of the little intervals to be the sizes the
lobes, but this map is
 only an operad map up to homotopy due to the
different scalings. To remedy the situation, we could augment the
spineless cacti operad to a larger homotopy equivalent model. In
this model the outside circle will be additional data. It will
have to be orientation preserving but not necessarily injective.
The parameterizations will be allowed to have stops at the
intersections points and the start (the global zero).

\subsection{The Gerstenhaber structure}
Due to the theorem of Cohen \cite{cohen, cohen2} identifying Gerstenhaber
algebras with algebras over the homology of the little discs
operad and the Theorem \ref{cactdisc} above, we know that algebras
over the homology of the operad of spineless cacti are
Gerstenhaber algebras.

\subsubsection{The explicit presentation of the operations in normalized
spineless cacti} 
\begin{figure}[ht!]
\cl{\epsfbox{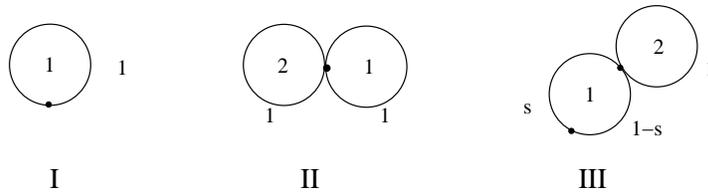}}
\caption{\label{ops}I: The identity\qua II: The product $\cdot$\qua
III: The operation $*$}
\end{figure}

Parallel to \cite{CS,V,KLP}, we can give explicit
generators for the operations on the chain level yielding the
Gerstenhaber structure on the homology spineless cacti (see Figure
\ref{ops}).

\begin{figure}[ht!]
\epsfxsize = \textwidth
\epsfbox{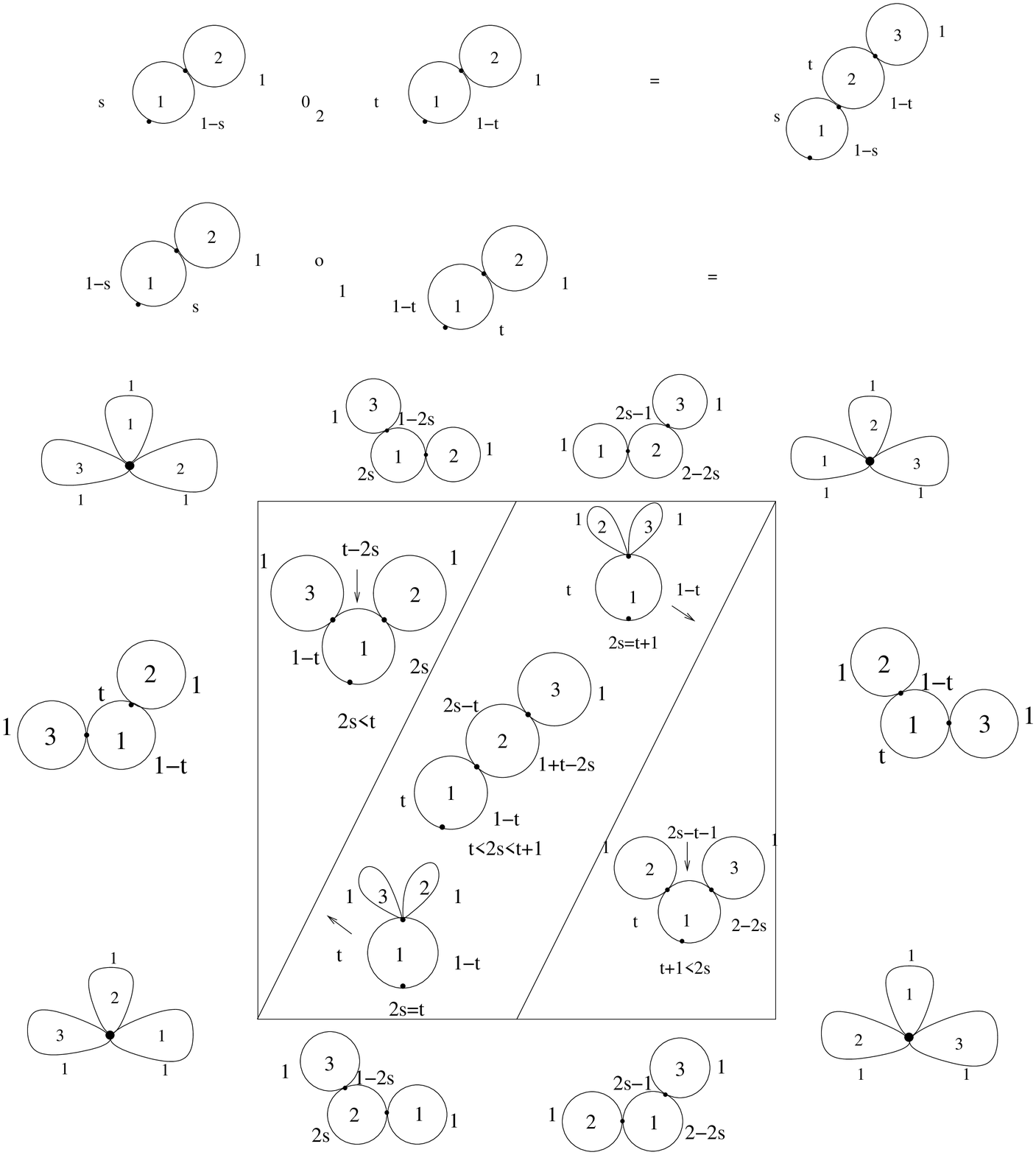}
\caption{\label{cactcircs}The associator in  normalized spineless cacti}
\end{figure}

We would like to emphasize that the product $\cdot$ is {\em
associative} on the nose already on the chain level. As usual, the
multiplication $*$ defines the bracket via the odd commutator.
\begin{equation}
\label{bracketdef}
\{a , b\}:= a*b - (-1)^{(|a|+1)(|b|+1)}b*a
\end{equation}
where we denoted the degree of $a$ and $b$ by $|a|$ and $|b|$. Its
iterations are given in Figure \ref{cactcircs} from which one can
also read off the associator (pre-Lie) relation which guarantees
the odd Jacobi identity.

Using the dual graph construction of Appendix B all the other
chain homotopies can be made explicit by translating them from
\cite{KLP} to normalized cacti.

 Just like Corollary \ref{eone} one obtains:
\subsubsection{Corollary}
\label{etwo}
{\sl The operad of the decorated moduli space of bordered,
punctured surfaces with marked points on the boundary $\widetilde
M^s_{g,r}$ which is proper homotopy equivalent to
 $\Arc_{\#}$ contains an $E_2$ operad. Thus so does the operad
 $\Arc$.

 The same is true for the
 operad of the moduli spaces $M_{g,n}^n$ of genus $g$ curves with $n$
 punctures and $n$ tangent direction and its restriction to genus
 $0$.

 Finally the spaces $M^1_{g,n}$ of surfaces of genus $g$, $n$ marked points
and a tangent vector at the first marked point form a partial
operad which is an
 $E_2$ operad.}

\begin{proof}
The proof of all but the last statement is analogous to the proof
of Corollary \ref{eone}. For the last statement we need to use the
fact that there is a cell model for $M^1_{g,n}$ by ribbon graphs
with one marking (cf. \cite{hoch}) and that the rooted tree-like
ribbon graphs given by $\Cact$ form an operad structure which
gives partial operad on the whole space.
\end{proof}

\subsubsection{Remark} Using Cohen's Theorem
this means that on the cell level algebras over these operads will
be Gerstenhaber algebras up to homotopy and on the homology level
Gerstenhaber algebras. In particular the operads themselves
possess these properties.

The operad $M_{g,n}$ is not included in Corollary \ref{etwo},
since we do need the markings of a global root on the ribbon
graphs in order to define the gluing. For spineless cacti one such
marking is enough, however.

\subsubsection{Remark}
We wish to point out that the chain defining the product $*$ is
exactly the path $\a_1$ for $\Cact(2)$ and
 diagram for the associator Figure
\ref{cactcircs} coincides up to re-parametrization with the braid
relation for the paths $\a_1$ and $\a_2$ in $\Cact(3)$.

One can obtain all paths $\a$ and all braid relation by taking
quasi-operadic products with spineless corolla cacti whose radii
are all one. These can be viewed as a quasi-sub-operad of the
quasi-operad of normalized spineless cacti. In fact, the
base-points $c_n$ suffice for this, i.e\ the path $\a_i$ in
$\Cact(n)$ is $c_{n-2} \circ_i \a_1$ and likewise one obtains the
braid relation.

This nicely ties together the point of view of \cite{fiedorbar} and
\cite{cohen, cohen2} in relating the Gerstenhaber structure directly to
the braided structure.

\section{Examples and constructions}
\label{prods}

In this section, we collect constructions and results which we
will modify and use in section \ref{relationssection} to relate
our different varieties of cacti in terms of semi-direct and
bi-crossed products.

\subsection{Operads of spaces}
\label{spaceop} The following procedure is motivated by
topological spaces with Cartesian product, but actually works in
any strict symmetric monoidal category where the monoidal product
is a product (i.e. we have projection maps).

Let $X$ be a topological space, then we can form the iterated
Cartesian product $X\times \dots \times X$. We simply denote the
$n$-fold product by $X(n)$. This space has an action of $\Sn$ by
permutation of the factors. We denote the corresponding morphism
also by elements of $\Sn$. Given a subset $I \subset\{1,
\dots,n\}$ we denote the projection $\pi_I: X(n) \rightarrow X(I)=
\times_{i\in I} X$.
\begin{multline}
\bar\circ_i: X(n) \times X(m) \stackrel{\pi_{\{1,\dots,n\}
\setminus i}}
{\rightarrow} X(n-1) \times X(m) \rightarrow X(m+n-1)\\
\stackrel{\s}{\rightarrow} X(m+n-1)
\end{multline}
where $\s\in {\Bbb S}_{m+n-1}$ is the permutation that shuffles
the last $m$ factors into the place $i$. I.e.\
\begin{equation}
((x_1, \dots,x_n)\bar\circ_i (x'_1,\dots,x'_m))= (x_1,\dots
x_{i-1},x'_1,\dots , x'_m,x_{i+1},\dots, x_n).
\end{equation}

\subsubsection{The cyclic version} Using $X((n))= X(n+1)$ with the
$\Snn$ action and the gluing above, one obtains a cyclic version
of the construction.

\subsection{Operads built on monoids}
Let $S$ be a monoid with associative multiplication $\mu: S \times
S \rightarrow S$. For simplicity we will denote this
multiplication just by juxtaposition: $s,s' \in S; ss':=
\mu(s,s')$. We will take $S$ to be an object in a strict symmetric
monoidal category.

We set $S(n):= S^{\times n}$ and endow it with the permutation
action.

\subsubsection{An operad defined by a monoid}
We consider the following maps:
\begin{eqnarray}
\circ_i: S(n)\times S(m)&\rightarrow&S(n+m-1)\nn\\
((s_1, \dots,s_n),(s'_1,\dots,s'_m))&\mapsto& (s_1,\dots
s_{i-1},s_i s'_1,\dots ,s_i s'_m,s_{i+1},\dots, s_n)\nn
\end{eqnarray}
It is straightforward to check that these maps define an operad in
the same category as $S$. This operad is unital if $S$ is unital.

\subsubsection{Examples}$\phantom{99}$
\begin{itemize}
\item[1)] One standard example is that of a Lie group in
topological spaces.

\item[2)] Another nice example is that of a field $k$. Then
$k(n)=k^n$ and the gluing is plugging in vectors into vectors
scaled by scalar multiplication. There are ${\mathbb Z}$-graded
and super versions of this given by including the standard
supersign for the permutation action where in the ${\mathbb
Z}$-graded version one uses the induced $\ZZ$-grading.

\item[3)] The example $S^1$ is particularly nice. In this case
(see e.g.\ \cite{KLP}), we can see that for the monoid $S^1$ the
homology operad $H_*(S^1(n),\ZZ)$ of its induced operad $S^1(n)$
is isomorphic, though not naturally, to the operad built on $\ZZ$
and that for a field $k$, $H_*(S^1(n),k)$ is isomorphic to the
direct product of the operads $Comm$ and the operad built on the
monoid $\ZZ$.
\end{itemize}

\subsubsection{Remark} There are several other natural versions of
operads and cyclic operads which can be defined analogously to the
operads built on circles which are presented in \cite{KLP}.

\subsection{Semi-direct products with monoids}
We now turn to the situation where the monoid $S$ acts on all the
components of an operad.

I.e.\ Let $Op(n)$ be an operad in a symmetric monoidal category
$\mathcal{C}$ and let $S$ be a monoid in the same category such
that $S$ acts on $Op(n)$.
$$
S \times Op(n) \stackrel{\rho}{\rightarrow} Op(n)
$$
s.t.\ the following diagrams are commutative
$$
\begin{CD}
S\times S \times Op(n)@>\mu\times id>> S\times Op(n)\\
@A id\times \rho AA @VV \rho V\\
S\times Op(n) @>\rho>> Op(n)\\
\end{CD}
$$
$$
\begin{CD}
S\times Op(m) \times Op(n)@>id\times \circ_i>> S\times Op(m+n-1)@>\rho>>Op(m+n-1)\\
@V \Delta\times id \times id VV && @A\circ_iAA\\
S{\times} S {\times} Op(m) {\times} Op(n) @>\sigma_{23}>> S {\times} Op(m){\times} S {\times} Op(n)
@>\rho{\times} \rho>> Op(m){\times} Op(n)\\
\end{CD}
$$
Consider the action
\begin{eqnarray}
\rho_i: S(n) \times Op(m)&\rightarrow& Op(m)\nn\\
((s_1, \dots, s_n),o) &\mapsto&\rho(s_i)(o)
\end{eqnarray}
and the twisted multiplications
\begin{eqnarray}
\circ^s_i: Op(n) \times Op(m) &\rightarrow& Op(n+m-1)\\
(o,  o') &=& (O \circ_i \rho(s,O'))
\end{eqnarray}
It is straightforward to check that the action $\rho$ satisfies
the conditions of Lemma \ref{actionlemma}.

\subsubsection{Definition}
We define the semi-direct product $Op\rtimes S$ of an operad $Op$ with
a monoid $S$ in the same category to be given by the spaces
$$
(Op \rtimes S)(n) := Op(n) \times S(n)
$$
with diagonal $\Sn$ action and the compositions
\begin{eqnarray}
\label{semi} \circ^{\rtimes}_i:(Op \rtimes S)(n) \times (Op
\rtimes S)(m)&\rightarrow& (Op
\rtimes S)(n+m-1)\nn\\
((O,s),(O',s'))&:=& (O \circ^{\pi_i(s)}_i O',s\circ_i s')
\end{eqnarray}
where $\pi_i$ is the projection to the i-th component.

\subsubsection{Proposition}{\sl Given an action $\rho$ and the 
multiplications $\circ i^{\rtimes}$ defining the quasi-operad
structure of the semi-direct product as above, the semi-direct product
quasi-operad is an operad.}

\subsubsection{Example} The operad of framed
little discs is the semi-direct product of the little discs
operad with the operad based on the monoid given by the
 circle group $S^1$ \cite{SW}.

\subsubsection{The action} If we are in a category that
satisfies the conditions of \ref{spaceop}, we can break down the
operad structure based on a monoid into two parts. The first is
the structure of operads of spaces and the second is the diagonal
action.

More precisely let $\Delta: S \rightarrow S(n)$ be the diagonal
and $\mu: S\times S \rightarrow S$ be the multiplication:
$$
\rho_{\Delta}: S \times S(n) \stackrel{\Delta}{\rightarrow} S(n)
\times S^n  \stackrel{\mu^n}{\rightarrow} S(n)
$$
Here we denote by $\mu^n$ the diagonal multiplication:
$$
\mu^n((s_1, \dots, s_n),(s'_1, \dots, s'_n))=(s_1s'_1, \dots,
s_ns'_n)
$$
\begin{multline}
\circ_i: S(n) \times S(m) \stackrel{(id \times \pi_i)(\Delta)
\times id} {\longrightarrow}
S(n) \times S \times S(m)\\
\stackrel{id \times \rho_{\Delta}} {\longrightarrow} S(n) \times
S(m)
 \stackrel {\bar\circ_i}{\rightarrow}
S(n+m-1)
\end{multline}
where $\bar\circ_i$ is the operation of the operad of spaces and
$$
(id \times \pi_i)(\Delta): S(n) \stackrel{\Delta}{\rightarrow}
S(n) \times S(n) \stackrel{id \times \pi_i}{\rightarrow} S(n)
\times S.
$$

\section{The relations of the cacti operads}
\label{relationssection}

\subsection{The relation between normalized (spineless) cacti and\newline
(spineless) cacti}

\subsubsection{The scaling operad}
We define the scaling operad  ${\mathcal R}_{> 0}$ to be given by the
spaces ${\mathcal R}_{> 0}(n):=  {\Bbb R}_{> 0}^{n}$ with the
permutation action by $\Sn$ and the following products
$$
(r_1,\dots,r_n) \circ_i (r'_1,\dots,r'_m)
= (r_1, \dots r_{i-1}, \frac{r_i}{R} r'_1,
\dots, \frac{r_i}{R} r'_m, r_{i+1},\dots r_{n})
$$
where $R=\sum_{k=1}^m r'_k$. It is straightforward to check that
this indeed defines an operad.

\subsection{The perturbed compositions}
We define the perturbed compositions
\begin{equation}
\label{percactmult} \circ_i^{{\mathcal R}_{>0}}: \Cacti^1(n)
\times {\mathcal R}_{> 0}(m) \times  \Cacti^1(m) \rightarrow
\Cacti^1(n+m-1)
\end{equation}
via the following procedure: Given $(c,\vec{r}',c')$ we first scale
$c'$ according to $\vec{r}'$, i.e.\ scale the $j$-th lobe of $c'$ by
the $j$-th entry $r_i$ of $\vec{r}$ for all lobes. Then we scale the
$i$-th lobe of the cactus $c$ by $R=\sum_j r_j$ and glue in the
scaled cactus. Finally we scale all the lobes of the composed cactus back to one.

We also use the analogous perturbed compositions for $\Cact^1$.

\subsubsection{The perturbed multiplications in terms of an action}
We can also describe, slightly more technically, the
above compositions in the following form.
Fix an element $\vec{r} :=(r_1,\dots,r_n) \in {\Bbb R}^n_{> 0}$
and set $R=\sum_i r_i$ and a normalized cactus $c$ with $n$ lobes. Denote by
$\vec{r}(c)$ the cactus where each lobe has been scaled according to $\vec{r}$,
i.e.\ the $j$-th lobe by the $j$-th entry of $\vec{r}$. Now consider the chord diagram
of the cactus $\vec{r}(c)$. It defines an action on $S^1$
via
\begin{equation}
\rho: S^1 \stackrel {rep^1_R}{\longrightarrow}S^1_R
\stackrel{cont_{\vec{r}(c)}}{\longrightarrow} S^1_n
\stackrel {rep^n_1}{\longrightarrow}  S^1
\end{equation}
Where $cont_{\vec{r}}$ acts on $S^1_R$ in the following way.
Identify the pointed $S^1_R$ with the pointed outside circle of the chord diagram of $\vec{r}(c)$.
Now contract the arcs belonging to the $i$-th lobe  homogeneously
with a scaling factor $\frac{1}{r_i}$.

We think of the $i$-lobe of a normalized (spineless) cactus as an
$S^1$ with base point given by the local zero together with
additional marked points;  the additional marked points are the
intersection points. Using the map above on the $i$-lobe we thus
obtain maps:
\begin{eqnarray}
\rho_i: \Cact^1(n) \times {\mathcal R}_{> 0}(m) \times \Cact^1(m)
&\rightarrow& \Cact^1(n) \nn\\
\rho_i: \Cacti^1(n) \times {\mathcal R}_{> 0}(m)\times \Cacti^1(m)
&\rightarrow& \Cacti^1(n)
\end{eqnarray}
What this action effectively does is move the
lobes and if applicable the root of the cactus $c$ which are
attached to the i-th lobe according to the
cactus $\vec{r}(c')$ in a manner that depends continuously on
$\vec{r}$ and $c'$.

With this action we can write the perturbed multiplication as:
\begin{multline}
\circ_i^{{\mathcal R}_{>0}}: \Cacti^1(n) \times
{\mathcal R}_{> 0}(m) \times  \Cacti^1(m)\\
\stackrel{id\times id \times \Delta}{\longrightarrow} \Cacti^1(n)
\times {\mathcal R}_{>
0}(m) \times  \Cacti^1(m)\times  \Cacti^1(m)\\
 \stackrel{\rho_i \times id}{\longrightarrow}
\Cacti^1(n) \times \Cacti^1(m) \stackrel{\circ_i}{\longrightarrow}
\Cacti^1(n+m-1)
\end{multline}

\subsubsection{Theorem}
\label{semiprodthm}{\sl The operad of spineless cacti is isomorphic to
the operad given by the semi-direct product of their normalized
version with the scaling operad. The latter is homotopy equivalent
(through quasi-operads)  to the direct product as a quasi-operad.
The direct product is in turn equivalent as a quasi-operad to
$\Cact^1$. The same statements hold true for cacti.
\begin{eqnarray}
\Cact &\cong&   {\mathcal R}_{>0}\ltimes \Cact^1
\sim  \Cact^1 \times {\mathcal R}_{>0}\simeq \Cact^1\nn\\
\Cacti &\cong&  {\mathcal R}_{>0} \ltimes \Cacti^1 \sim \Cacti^1
\times {\mathcal R}_{>0}\simeq \Cacti^1
\end{eqnarray}
here the semi-direct product compositions are given by:
\begin{equation}
(\vec{r},c)\circ_i(\vec{r}',c') = ( \vec{r}\circ_i \vec{r'}, c
\circ_i^{\vec{r'}} c')
\end{equation}}

\begin{proof} By definition the
space $\Cact(n)=\Cact^1(n)\times \mathbb{R}^n_{>0}=
\Cact^1(n)\times \mathcal{R}_{>0}(n)$. To establish that the
operad structure of $\Cact$ is that of a semi-direct product,
first notice that the behavior of the radii under gluing is given
by the scaling operad. Second we notice that the global incidence
conditions, i.e.\ the positions of the intersection points on the
outer circle, are shifted under the scaling in a way that is
compensated by perturbed multiplication. This means that when
gluing, we do not use the outside circle of the normalized cactus,
but the outside circle of the original cactus which is recovered
form the outside circle of the normalized cactus and the radii by
the action of the scaling operad.

Secondly, the perturbed multiplications are homotopic to the
unperturbed multiplications: a homotopy is for instance given by a
the choice of a  paths $\vec{r}_t$ given by the line segment from
$\vec{r}$ to $(1,\dots,1)$ and the following homotopy for the
quasi-operadic compositions
$$
(\vec{r},c)\circ_{i,t}(\vec{r}',c') = ( \vec{r}\circ_i \vec{r'}, c
\circ_i^{\vec{r}_t'} c')
$$
and the composition $\circ_i^{\vec{r}_t'}$ uses the action chord diagram
of $\vec{r}_t'(C')$, i.e.\ the cactus $C'$ scaled by $\vec{r}'_t$, for the gluing.

Lastly, the $\mathcal{R}(n)=\mathbb{R}^n_{>0}$ are contractible
and the contraction induces an equivalence of quasi-operads.

The analogous arguments hold for $\Cacti$.
\end{proof}

\subsubsection{Corollary}
{\sl Cacti without spines are homotopy equivalent to normalized cacti
without spines as quasi-operads. The  quasi-operad of homology
normalized spineless cacti is an operad which is isomorphic to the
homology of the spineless cacti operad.

Also, cacti and normalized cacti are homotopy equivalent as
quasi-operads and the homology quasi-operad of normalized cacti is
an operad which is isomorphic to the homology of the cacti operad.}

\begin{proof}
We first remark that there is a homotopy of quasi-operads to the direct product,
then there is a homotopy of quasi-operads from the scaling operad to the operad of
spaces. Finally there is a homotopy between the operad of spaces built on $\mathbb{R}_{>0}$
and a point in each degree.
\end{proof}

From the previous analysis, we obtain:
\subsubsection{Corollary}
{\sl The quasi-operads of normalized cacti and normalized\break spineless cacti
are both homotopy associative.}

\subsubsection{Remark} It is shown in \cite{del} that the quasi-operad structure
induced on the cellular chains $CC_*(\Cact^1)$ is actually an operad structure. This can
be seen from the explicit description of the semi-direct product above.
It follows from the Theorem above and Theorem \ref{cactdisc} that $CC_*(\Cact^1)$
is a cell model for the little discs operad.

The fact that the cells are indexed by trees then yields a quick proof of Deligne's conjecture \cite{del}
which states that the Hochschild cochains of an associative algebra are an algebra over
a cell model operad of the little discs.

\subsection{The relation between spineless cacti and cacti}
We would now like to specialize the monoid of \S \ref{prods}
to $S= S^1$. We already showed
that the normalized and non-normalized versions of the different
species of cacti are homotopic and moreover that they are related by
taking the direct product with the scaling operad.
Below we will see that cacti with spines are a bi-crossed product
of the cacti without spines and the operad built on $S^1$. Furthermore
we show that this bi-crossed product is homotopic to the semi-direct
product. Thus we see that the relation of cacti with and without
spines is analogous to the relation of framed little discs and little discs.

\subsubsection{The action of $S^1$ and the twisted gluing}
There is an action of $S^1$ on $\Cact(n)$ given by rotating the
base point {\em clockwise} around the perimeter. We denote this
action by
$$\rho^{S^1}: S^1 \times \Cact(n) \rightarrow \Cact(n).$$
With this action we can define the twisted gluing:
\begin{eqnarray}
\label{circtheta}
\circ_i^{S^1}:\Cact(n) \times S^1(n) \times \Cact(m) &\rightarrow& \Cact(n+m-1)\nn\\
(C,\theta,C')&\mapsto& C \circ \rho^{S^1}(\theta_i,C') =: C \circ_i^{\theta_i}C'
\end{eqnarray}

\subsubsection{The homotopy diagonal defined by a spineless cactus}
\begin{figure}[ht!]
\cl{\epsfxsize = .6\textwidth
\epsfbox{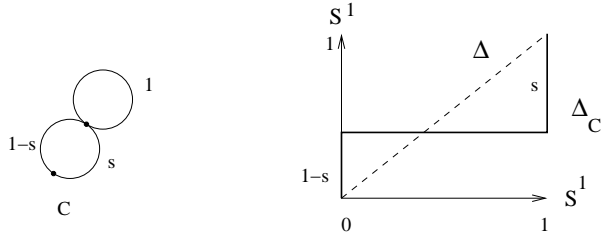}}
\caption{\label{cdiagfig}
The image of $\Delta_C$ for a two component cactus}
\end{figure}

Given a cactus  without spines $C\in \Cact(n)$ the orientation
reversed perimeter (i.e.\ going around the outer circle {\em
clockwise}) gives a map
\begin{equation}
\Delta_C: S^1 \rightarrow (S^1)^n.
\end{equation}
As one goes around the perimeter the map goes around each circle
once and thus the map $\Delta_C$ is homotopic to the diagonal
\begin{equation}
\label{homotop}
\Delta_C (S^1) \sim \Delta(S^1).
\end{equation}
A picture of the image of $\Delta_C$ for a two component cactus is
depicted in Figure \ref{cdiagfig}.

\subsubsection{The action based on a cactus}
We can use the map $\Delta_C$ to give an action of $S^1$ and
$(S^1)^{\times n}$.
\begin{equation}
\rho^C: S^1 \times(S^1)^{\times n}\stackrel{\Delta_C}
{\rightarrow} (S^1)^{\times n} \times (S^1)^{\times n}
\stackrel{\mu^n}{\rightarrow}(S^1)^{\times n}
\end{equation}
And furthermore using concatenations with projections we can
define maps
\begin{multline}
\label{perturbdef}
\circ^C:(S^1)^{\times n} \times (S^1)^{\times m} \stackrel{(id
\times \pi_i)(\Delta) \times id} {\longrightarrow} (S^1)^{\times n} \times
S^1\times (S^1)^{\times m}\\
\stackrel{id \times \rho^C}{\longrightarrow}
(S^1)^{\times n} \times (S^1)^{\times m}
\stackrel{\bar\circ_i}{\longrightarrow}(S^1)^{\times n+m-1}
\end{multline}

\subsubsection{Theorem}
\label{bicrossed}{\sl The (quasi-)operad of (normalized) cacti is
isomorphic to the bi-crossed product of the operad of (normalized)
spineless cacti with the operad ${\mathcal S}^1$ based on $S^1$
with respect to the compositions of (\ref{circtheta}) and
(\ref{perturbdef}) and furthermore this bi-crossed product is
homotopy equivalent as a  quasi-operad to the semi-direct product
of the operad of cacti without spines with the circle group $S^1$
considered as a monoid.
\begin{equation}
\Cacti \cong \Cact \bowtie {\mathcal S}^1 \sim \Cact \rtimes
{\mathcal S}^1
\end{equation}
with respect to the operations of (\ref{circtheta}) and (\ref{perturbdef}).
Similarly
\begin{equation}
\Cacti^1 \cong \Cact^1 \bowtie {\mathcal S}^1 \sim \Cact \rtimes
{\mathcal S}^1
\end{equation}}

\begin{proof} As spaces the operad of cacti is the direct product
of spineless cacti and the operad built on the monoid $S^1$
$$
\Cacti(n) = \Cact(n) \times S^1(n)
$$
where for the identification we use intersection points and the
global zero to define the parameterizations of the $S^1$s
constituting the cactus. Then the  local zeros are specified by
their coordinate, i.e.\ all zeros are fixed by a point on
$S^1(n)$.

The product on the cacti operad in this identification is given by
\begin{equation}
(C,\theta) \circ_i (C',\theta') = (C\circ_i^{\theta_i} C',
\theta\circ_{i}^{C'}\theta')
\end{equation}
where we used the operations of (\ref{circtheta}) and (\ref{perturbdef}).

This comes from the observation that gluing the global zero to the
local zero is the same as first using the $S^1$ action that
rotates the global zero around the perimeter in a clockwise
fashion by the amount which is given by the coordinate of the
local zero on the cactus to be glued in and then using the
standard gluing. During the movement of the global zero the
coordinates of the local zeros change according to the map
$\rho_C$. This proves the first assertion.

Now, in the semi-direct product the multiplication is given by:
\begin{equation}
(C,\theta) \circ_i (C',\theta') = (C\circ_i^{\theta_i} C',
\theta\circ_{i}\theta')
\end{equation}
Now by (\ref{homotop}) and (\ref{perturbdef}) $\circ_i^{C} \sim
\circ_i$ are homotopic and $\circ_i^{C}$  depends continuously on
$C$. Choosing the following simultaneous homotopies of the
$\Delta_C$ with $\Delta$, we obtain a homotopy equivalence of
quasi-operads. The homotopy is given by straightening out $\Delta
C$. For $\Delta_C(\theta)=(\theta_1,\dots, \theta_n)$ notice that
$\theta= \frac{1}{n}\sum_{i=1}^n\theta_i$; we set
$$h_t(\theta)=((1-t)\theta_1+t\theta,
\dots,(1-t)\theta_n+t\theta)  \quad t\in [0,1].$$
This proves the
second statement. On the level of homology we therefore obtain an
isomorphism of operads, so that the map of $\Cacti$ to the
semi-direct product is a quasi-isomorphism and hence an
equivalence. To obtain the version of the theorem in the
normalized situation, we remark that the bi-crossed structure does
not depend on the size of the lobes.
\end{proof}

\subsubsection{Corollary}
\label{equicor}{\sl The operads of cacti and the semi-direct product
of spineless cacti with the operad $\mathcal{S}^1$ are  homotopy
equivalent as quasi-operads.}

\medskip

To connect these results to the literature, we recall a theorem of
Voronov announced in \cite{V,Vman}:
\subsubsection{Theorem}
\label{vorthm}{\sl The cacti operad is  equivalent  to the framed
little discs operad.}

\subsubsection{Proposition}
{\sl Theorem \ref{cactdisc} and the Corollary \ref{equicor} above in
conjunction with the equivariant recognition principle of
\cite{SW} imply Theorem \ref{vorthm}.}

\begin{proof}
For this first remark that by what we have proven above the spaces
$\Cacti$ already have the right weak homotopy type, they are
$K(PRBr,1)$. Here $PRBr$ is the pure ribbon braid group.  Next,
it is clear from the decomposition of the spaces
$\Cacti(n)=\Cact(n)\times (S^1)^n$, the results on the quotient
$\Cact(n)/\Sn$ and the fact that $\Sn$ acts by permutation on the
factors $S^1$ that $\Cacti(n)/\Sn$ is a $K(RBr,1)$.  Here $RBr$ is
the ribbon braid group. The construction of the ribbon braid
operad on the covers is done by using a lifting of the appropriate
paths to define the operad compositions and an action of the pure
ribbon braid group which covers the symmetric group action. Here
SCC again act as base-points. This will endow the universal covers of $\Cacti$
with a ribbon braid group
operad structure. This can be checked directly. It also follows from the
existence of the homotopies to the semi-direct product, since the
same holds true for the semi-direct product of the little discs
with $S^1$ and the little discs are equivalent to $\Cact$. Also,
the covers are contractible and  the ribbon braid group action is
free so that the covers form an $R_{\infty}$ operad and hence by
\cite{SW} $\Cacti$ is weakly homotopy equivalent to the framed
little discs. The fact that $\Cacti$ are homotopy equivalent to a
CW complex then allows one to upgrade from the weak equivalence to
an honest one.
\end{proof}

Vice-versa using the recognition principle of \cite{SW} and starting with
characterization of $\Cacti$ in terms of $\Cact$ above:

\subsubsection{Proposition}{\sl The Theorem \ref{vorthm} together with the
Theorem \ref{bicrossed} imply Theorem \ref{cactdisc}, namely that
the operad of (normalized) cacti without spines is
 equivalent to the little discs operad.}

\begin{proof}[Sketch of a proof] The main idea is
that the operad of little discs is embedded into the framed little
discs by setting the coordinates of the factors $S^1$ equal to
zero in the semi-direct product. The argument is as follows. Following
\cite{SW} one can show that if the universal cover of a
semi-direct product $fD=D\rtimes S^1$ is an $R_{\infty}$ operad
then the universal cover of the suboperad $D(n)\times 0^n$ is a
$B_{\infty}$ operad. Since by the above theorem $\Cacti$ are
weakly homotopy equivalent to the framed little discs and also the
semi-direct product of $\Cact$ with $S^1$ is equivalent with
$\Cacti$ it follows by \cite{SW} that the universal cover of the
semi-direct product is an $R_{\infty}$ operad, and hence the
universal cover of the natural inclusion of $\Cact$ into the
semi-direct product is a $B_{\infty}$ operad which implies that
$\Cact$ is an $E_2$ operad.
\end{proof}

\subsubsection{Remark}
It is shown in \cite{cyclic} that the quasi-operad
on the cellular chains $CC_*(\Cacti^1)$ is an operad. Due to the above result,
this operad agrees on the homology level with the semi-direct product of $\Cact$ with $S^1$ which
in turn agrees with the semi-direct product of the little discs operad with $S^1$. Thus on
the homology level we obtain an operad isomorphic to the homology of the framed little
discs operad and thus $CC_*(\Cacti^1)$ gives a chain level operad model for the framed little discs.

Again, the combinatorial description of the cells allows one to prove a theorem in the spirit of
Deligne's conjecture. Namely \cite{cyclic}, the Hochschild cochains of a Frobenius algebra are
an algebra over a cell level model of the framed little discs operad.

\small
\renewcommand{\theequation}{A-\arabic{equation}}
\renewcommand{\thesection}{A}
\setcounter{equation}{0}  
\setcounter{subsection}{0}
\section*{Appendix A: Graphs}
\addcontentsline{toc}{section}{Appendix A: Graphs}
\label{cellapp}

In this appendix, we formally introduce the graphs and the
operations on graphs used in our analysis of cacti. This language
is useful for describing the $\Arc$ operad which is done in
Appendix B. We also give a supplemental result characterizing
cacti as a certain type of ribbon graph. A cactus is a marked
treelike ribbon graph with a metric.

\subsection{Graphs} A graph $\Gamma$ is a tuple
$(V_{\Gamma},F_{\Gamma}, \imath_{\Gamma}: F_{\Gamma}\rightarrow
F_{\Gamma},\del_{\Gamma}:F_{\Gamma} \rightarrow V_{\Gamma})$ where
$\imath_{\Gamma}$ is an involution $i_{\Gamma}^2=id$ without fixed
points. We call $V_{\Gamma}$ the vertices of $\Gamma$ and
$F_{\Gamma}$ the flags of $\Gamma$. The edges $E_{\Gamma}$ of
$\Gamma$ are the orbits of the flags under the involution. A
directed edge is an edge together with an order of the two flags
which define it.

In case there is no risk of confusion we will drop the subscripts
$\Gamma$.

Notice that $f\mapsto (f,\imath(f))$ gives a bijection between
flags and directed edges.

We also call $F_{\gamma}(v):=\del^{-1}(v)\subset F_{\Gamma}$ the
set of flags of the vertex $v$ and call $|F_{\Gamma}(v)|$ the
valency of $v$ and denote it by $|v|$.

The geometric realization of a graph is given by considering  each
flag as a half-edge and gluing the half-edges together using
$\imath$. This yields a one-dimensional CW complex whose
realization we call the realization of the graph.

\subsection{Trees}
 A graph is connected if its realization is. A graph is
a tree if it is connected and its realization is contractible.

A rooted tree is a pair $(\t,v_0)$ where $\t$ is a tree and
$v_0\in V_{\t}$ is a distinguished vertex.

In a rooted tree there is a natural orientation for edges, in
which the edge points toward the root. That is we say $(f,\imath
(f))$ is naturally oriented if $\del(\imath(f))$ is on the unique
shortest path from $\del(f)$ to the root.

A bi-colored or black and white tree is a tree $\t$ together with
a map $\color:V\rightarrow \mathbb{Z}/2\mathbb{Z}$. Such a tree is
called bipartite if for all $f\in
F_{\t}:\color(\del(f))+\color(\del(\imath(f)))=1$, that is edges
are only between black and white vertices. We call the set
$E_w:=\color^{-1}(1)$ the white vertices. If $(f,\imath (f))$ is a
naturally oriented edge, we call the edge white if
$\del(\imath(f))\in E_w$.

\subsection{Planar trees and ribbon graphs}

A ribbon graph is a connected graph whose vertices are of valency
at least two together with a cyclic order of the set of flags of
the vertex $v$ for every vertex $v$.

A tree with a cyclic order of the flags at each vertex is called
planar.

A  graph with a cyclic order of the flags at each vertex gives
rise to bijections $N_v:F_v\rightarrow F_v$ where $N_v(f)$ is the
next flag in the cyclic order and since $F=\amalg F_v$ to a map
$N:F\rightarrow F$.

The orbits of the map $N \circ \imath$ are called the cycles or
the boundaries of the graph. These sets have the induced cyclic
order.

Notice that each boundary can be seen as a cyclic sequence of
directed edges. The directions are as follows. Start with any flag
$f$ in the orbit. In the geometric realization go along this
half-edge starting from the vertex $\del(f)$, continue along the
second half-edge $\imath(f)$ until you reach the vertex
$\del(\imath(f))$ then continue starting along the flag
$N(\imath(f))$ and repeat.

A planar tree has only one cycle $c_0$.

A planted planar tree is a rooted planar tree $(\t,v_0)$ together
with a linear order of the set of flags at $v_0$. Such a tree has
a linear order of all flags as follows, let $f$ be the smallest
element of $\del^{-1}(v_0)$, then every flag
appears in $c_0$ and defining the flag $f$ to be the smallest gives a
linear order on the set of all flags. This linear order induces a linear order on
all oriented edges and on all oriented edges, by restricting to the
edges in the orientation opposite the natural orientation i.e.\ pointing away from the root.

The genus $g(\Gamma)$ of a ribbon graph $\Gamma$ is given by
$$2g(\Gamma)+2=|V_\Gamma|-|E_{\Gamma}|+\#cycles$$
The surface $\Sigma(\Gamma)$ of a ribbon graph $\Gamma$ is the
surface obtained from the realization of $\Gamma$ by thickening
the edges to ribbons. I.e.\ replace each 0-simplex $v$ by a closed
oriented disc $D(v)$ and each 1-simplex $e$ by $e\times I$
oriented in the standard fashion. Now glue the boundaries of
$e\times [-1,1]$ to the appropriate discs in their cyclic order
according to the orientations. Notice that the genus of
$\Sigma(\Gamma)$ is $g(\Gamma)$ and that $\Gamma$ is naturally embedded as
the spine of this surface.

\subsection{Treelike and marked ribbon graphs}

A ribbon graph together with a distinguished cycle $c_0$ is called
treelike if the graph is of genus $0$  such that for all other
cycles $c_i\neq c_0$ if $f\in c_i$ then $\imath(f)\in c_0$.
In other words each edge is traversed by the cycle $c_0$.
Therefore there is a cyclic order on all (non-directed) edges,
namely the cyclic order of $c_0$.

A marked ribbon graph is a
ribbon graph together with a map $\mk:\{cycles\} \rightarrow
F_{\Gamma}$ satisfying the conditions:
\begin{itemize}
\item[i)] For every cycle $c$ the directed edge $\mk(c)$ belongs
to the cycle.

\item[ii)] All vertices of valence two are in the image of $\mk$,
that is $\forall v,|v|=2$ implies  $v\in \del(Im(\mk))$.
\end{itemize}

Notice that on a marked treelike ribbon graph there is a linear
order on each of the cycles $c_i$. This order is defined by
upgrading the cyclic order to the linear order $\prec_i$ in which
$\mk(c_i)$ is the smallest element.

A marked treelike ribbon graph is called spineless, if:

\begin{itemize}
\item[i)] There is at most one vertex of valence $2$. If there is
such a vertex $v_0$ then $\del(mk(c_0))=v_{0}$.

\item[ii)] The induced linear orders on the $c_i$ are compatible
with that of $c_0$, i.e.\ $f\prec_i f'$ if and only if
$\imath(f')\prec_0 \imath(f)$.
\end{itemize}

A metric $w_{\Gamma}$ for a graph is a map $E_v\rightarrow
\mathbb{R}_{>0}$.

\subsection{Graphs with a metric}
A projective metric for a graph is a class of metrics equivalent
under re-scaling, i.e.\ $w\sim w'$ if $\exists \lambda \in
\mathbb{R}_{>0}\forall e\in E: w(e)=\lambda w'(e)$.

The length of a cycle is the sum of the length of its edges
$length(c)=\sum_{f\in c} w(\{f,\imath(f)\})$.

A metric for a treelike ribbon graph is called normalized if the
length of each non-distinguished cycle is $1$.

A projective metric for a treelike ribbon graph is called
normalized if it has a normalized representative.

\subsection{Marked ribbon graphs  with metric and maps of circles}
For a marked ribbon graph with a metric, let $c_i$ be its cycles,
let $|c_i|$ be their image in the realization and let $r_i$ be the
length of $c_i$. Then there are natural maps $S_{r_i}^1\rightarrow
|c_i|$ which map $S^1$ onto the cycle by starting at the vertex
$v_i:=\del(\mk(c_i))$ and go around the cycle mapping each point
$\theta\in S^1$ to the point at distance $\frac{\theta}{2\pi}r_i$
from $v_i$ along the cycle $c_i$.

\subsection{Contracting edges}
The contraction $(\bar V_{\Gamma}, \bar F_{\Gamma},\bar
\imath,\bar \del)$ of a graph
$(V_{\Gamma},F_{\Gamma},\imath,\del)$ with respect to an edge
$e=\{f,\imath(f)\}$ is defined as follows. Let $\sim$ be the
equivalence relation induced by $\del(f)\sim\del(\imath(f))$. Then
let $\bar V_{\Gamma}:=V_{\Gamma}/\sim$, $\bar
F_{\Gamma}=F_{\Gamma}\setminus\{f,\imath(f)\}$ and $\bar \imath:
\bar F_{\Gamma}\rightarrow \bar F_{\Gamma}, \bar\del: \bar
F_{\Gamma}\rightarrow \bar V_{\Gamma}$  be the induced maps.

For a marked ribbon graph, we define the marking of $(\bar
V_{\Gamma}, \bar F_{\Gamma},\bar \imath,\bar \del)$ to be
$\overline{\mk}(\bar c)=\overline{\mk(c)}$ if
$\mk(c)\notin\{f,\imath(f)\}$ and $\overline{\mk}(\bar
c)=\overline{N\circ \imath(\mk (c))}$ if
$\mk(c)\in\{f,\imath(f)\}$, viz.\ the image of the next flag in
the cycle.

\subsection{Labelling graphs}
By a labelling of the edges of a graph $\Gamma$ by a set $S$, we
simply mean a map $E_{\Gamma}\rightarrow S$.

A labelling of a ribbon graph $\Gamma$ by a set $S$ is a map
$\{$cycles of $\Gamma\}\rightarrow S$.

By a labelling of a black and white tree  by a set $S$ we mean  a
map $E_w\rightarrow S$.

\subsection{Cacti as ribbon graphs}

By considering its image a cactus is naturally a marked treelike
ribbon graph, with a metric. Vice-versa, given such a graph, one
obtains a cactus.

\subsubsection{Proposition}{\sl A cactus with $n$ lobes
is equivalent to an $\{0,1, \dots ,n\}$ labelled marked treelike
ribbon graph with a metric. I.e.\ The set  $\Cacti(n)$ is  in
bijection with the respective set of graphs. The conditions of
being normalized and/or spineless are  compatible with this
bijection.}

\begin{proof}
Given a cactus, its image is a ribbon graph. The vertices are the
marked points and the edges are the arcs. The flags being pairs
$(v,a)$ of a marked point $v$ and an arc $a$ ending at $v$ with
the obvious involution and map $\del$. At each vertex there is a
cyclic order induced from the one in the plane. This graph has
$n+1$ cycles. First each lobe marked by $i$ yields a cycle $c_i$
(going clockwise) and secondly the outside circle (going
counter-clockwise) is a cycle which we will call $c_0$.
Furthermore the outside circle goes through each edge and hence
the ribbon graph is treelike for the distinguished cycle $c_0$.
The labelling of the cycles is implicit in this description. The
marking is given by marking the unique flag $(v,f)$ on the cycle
$c_i$ for which $v$ is the local zero, for $i>0$ and the unique
flag $(v,f)$ on $c_0$ for which $v$ is the global zero. Finally
the metric is given by associating to an edge representing an arc,
the length of that arc.

In the reverse direction, given a treelike marked ribbon graph
with a metric $\Gamma$, consider the surface $\Sigma(\Gamma)$ of
$\Gamma$. This is a surface of genus $0$ with $n+1$ boundary
components. We embed this surface into the plane as a disc with
holes where the outside circle of the disc corresponds to the
cycle $c_0$. In this embedding, $\Gamma$ realized as the spine of
$\Sigma(\Gamma)$ is a cactus by considering the maps
$S^1_{r_i}\rightarrow c_i$ for $i>0$ as above.

It is clear that these maps induce bijections between the sets
$\Cacti(n)$ and the respective set of graphs.

To be normalized means that $r_i=1$ for $i>0$ in both cases.
Finally to be spineless in both cases also means the same. For
this notice that condition i) says that there is a global zero
which may or may not be an intersection point and that all local
zeros are at the intersection points or at the global zero.
Furthermore the root component of the cactus is fixed as the
unique cycle $c_i$ s.t.\ $\imath(\mk(c_0))\in c_i$. The condition
$b$ then insures that the local zero of the root component is the
global zero. The reverse direction is immediate.
\end{proof}

\renewcommand{\theequation}{B-\arabic{equation}}
\renewcommand{\thesection}{B}
\setcounter{equation}{0}  
\setcounter{subsection}{0}
\section*{Appendix B: The arc operad}  
\addcontentsline{toc}{section}{Appendix B: The arc operad}
\label{arcapp}

In this appendix, we would like to briefly recall some facts from
\cite{KLP} about the arc operad. We will make use of some
reformulations given in \cite{hoch} of the results of \cite{KLP}
which use the language of graphs, since this will simplify the
exposition.

\subsection{The arc operad}
We would like to recall some definitions of  \cite{KLP}. For this
we will fix an oriented surface $F_{g,r}^s$ of genus $g$ with $s$
punctures and $r$ boundary components which are labelled from $0$
to $r-1$, together with marked points on the boundary, one for
each boundary component. We call this data $F$ for short if no
confusion can arise.

We recall from \cite{KLP} that the space $A_{g,n}^s$ is the CW
complex whose cells are indexed by graphs on the surface $F$ up to
the action of the pure mapping class group $PMC$ which is the
group of elements of the mapping class group which fixes the
boundaries pointwise. A quick review in terms of graphs is as
follows.

\subsubsection{Embedded graphs}
By an embedding of a graph $\Gamma$ onto a surface $F$, we mean an
embedding  $i:|\Gamma|\rightarrow F$ with the conditions:

 \begin{itemize}
\item[i)] $\Gamma$ has at least one edge.

\item[ii)] The vertices map bijectively to the marked points on
the boundaries.

\item[iii)] No images of two edges are homotopic to each other.

\item [iv)]No image of an edge is homotopic to a part of the
boundary.
\end{itemize}

Two embeddings are equivalent if there is a homotopy of embeddings
of the above type  from one to the other. Note that such a
homotopy is necessarily constant on the vertices.

The images of the edges are called arcs. And the set of connected
components of $F\setminus j(\Gamma)$ are called complementary
regions.

Changing representatives in a class  yields a natural bijection of
the sets of arcs and connected components, we will therefore
associate arcs and connected components also with a class of
embeddings.

\subsubsection{Definition}
By a  graph on a surface we mean a triple $(F,\Gamma,[i])$ where
$[i]$ is an equivalence class of embeddings of $\Gamma$.

\subsubsection{A linear order on arcs}
Notice that due to the orientation of the surface the graph
inherits an induced linear order of all the flags at every vertex
$F(v)$ from the embedding. Furthermore there is even a linear
order on all flags by enumerating the flags first according to the
boundary components on which their vertex lies and then according
to the linear order at that vertex. This induces a linear order on
all edges by enumerating the edges by the first appearance of a
flag of that edge.

\subsubsection{The poset structure}
The set of such graphs on a fixed surface $F$ is a poset. The
partial order is given by calling $(F,\Gamma',[i'])\prec
(F,\Gamma,[i])$ if $\Gamma'$ is a subgraph of $\Gamma$ with the
same vertices and $[i']$ is the restriction of $[i]$ to $\Gamma'$.
In other words, the first graph is obtained from the second by
deleting some arcs.

 We
associate a simplex $\Delta(F,\Gamma,[i])$ to each such graph.
$\Delta$ is the simplex whose vertices are given by the set of
arcs/edges enumerated in their linear order. The face maps are
then given by deleting the respective arcs. This allows us to
construct a CW complex out of this poset.

\subsubsection{Definition}
Fix $F=F_{g,n}^s$. The space $A_{g,n}^{\prime s}$ is the space
obtained by gluing the simplices $\Delta(F,\Gamma',[i'])$ for all
graphs on the surface according to the face maps.

The pure mapping class group ($PMC$) which is the part of the
mapping class group that preserves the boundary pointwise
naturally acts on $A_{g,n}^{\prime s}$. It actually has finite
isotropy \cite{KLP}.

\subsubsection{Definition}
The space $A_{g,r}^s:= A_{g,r}^{\prime s}/PMC$.

\subsubsection{Elements of the $A_{g,r}^s$ as projectively weighted graphs}
The space $A_{g,n}^s$ is a CW complex whose cells are indexed by
graphs on the surface $F$ up to the action of the pure mapping
class group $PMC$. Moreover the cell for a given class of graphs
is actually a simplex whose vertices correspond to the arcs in the
order discussed above. The attaching maps are given by deleting
edges. Due to the action of $PMC$ some of the faces also might
become identified by the gluing.

Using barycentric coordinates the elements of $A_{g,n}^s$ are
graphs on the surface $F$ up to the action of the pure mapping
class group $PMC$ together with a map $w$ of the edges of the
graph $E_{\Gamma}$ to $\mathbb{R}_{>0}$ assigning a weight to each
edge s.t.\ the sum of all weights is 1.

Equivalently we can regard the map $w:E_{\Gamma}\rightarrow
\mathbb{R}_{>0}$ as an equivalence class under the equivalence
relation of, i.e.\ $w\sim w'$ if $\exists \lambda \in
\mathbb{R}_{>0}\forall e\in E_{\Gamma} \;  w(e)=\lambda w'(e)$.
That is $w$ is a projective metric. We call the $w(e)$ the
projective weights of the edge.
 In the limit, when the projective weight of an
edge goes to zero, the edge/arc is deleted.

\subsubsection{Example}
$A_{0,2}^0=S^1$. Up to PMC there is a unique graph with two edges.
This gives a one simplex. The two subgraphs lie in the same orbit
of PMC and thus the 0-simplices are identified to yield $S^1$. The
fundamental cycle is given by $\delta$ of Figure \ref{1ops}.

\subsubsection{Drawing pictures for arcs} There are several
pictures one can use to view elements in the arc operad. In order
to draw elements of the $\Arc$ operad it is useful, to expand the
marked point on the boundary to an interval, and let the arcs end
on this interval according to the linear order. Equivalently, one
can mark one point of the boundary and let the arcs end in their
linear order anywhere but on this point.

\subsubsection{Remark}
There is an operad structure on $\Arc$ which is obtained by gluing
the surfaces along the boundaries and splitting the arcs according
to their weights. We refer the reader to \cite{KLP} for the
details.

\subsubsection{Definition of suboperads and $\Darc$}$\phantom{99}$

Let $\Darc_{g,r}^s := \Arc_{g,r}^s\times \mathbb{R}_{>0}$.

Let  $\Arc_g^s(n)\subset A_{g,n+1}^s$ be the subspace of graphs on
the surface with the condition that $\del:F_{\Gamma}\rightarrow
V_{\Gamma}$ is surjective. This means that each boundary gets hit
by an arc.

Let  $\Arc^s_{\#g}(n)\subset \Arc_{g}^s(n)$ be the subspace of elements
whose complementary regions
are polygons or once punctured polygons.

Set $\Arc_{cp}(n):=\Arc_0^0(n)$.

Let $\Tree_{cp}(n) \subset \Arc_{cp}(n)$ be the subspace in which
all arcs run from $0$ to some boundary $i$ only.

Finally let $\Lintree_{cp}\subset \Tree_{cp}$ be the space in
which the linear order of  the arcs at the boundary $0$ is
anti-compatible with the linear order at each boundary, i.e. if
$\prec_i$ is the linear order at $i$ then if $f\prec_i f'$
$\imath(f')\prec_0 \imath(f)$.

\subsubsection{Remarks}$\phantom{99}$
\begin{itemize}

\item[1)] The elements of $\Darc$ are graphs on surfaces with a
metric, i.e.\ a function  $w:E_{\Gamma}\rightarrow
\mathbb{R}_{>0}$. And $\Darc$ is an operad equivalent to $\Arc$
\cite{KLP}.

\item[2)] The subspaces above are actually suboperads \cite{KLP}.

\item[3)] Any suboperad $\mathcal{S}$ of the list above defines a
suboperad $\mathcal{DS}:=\mathcal{S}\times \mathbb{R}_{>0}$ of
$\Darc$ which is equivalent to $\mathcal{S}$.

\item[4)] One can also reverse the orientation at zero. This is in
line with the usual cobordism point of view used in \cite{KLP}. In
this case the condition for $\Lintree$ is the compatibility of the
orders.

\end{itemize}

In \cite{KLP} we defined a map called $\mathcal{L}oop$ which is
the suitable notion of a dual graph for a graph on a surface. This
map is uses an interpretation of the graph as a partially measured
foliation. If one restricts to the subspace $\Arcn$ though, this
map has a simpler purely combinatorial description. This
description will be enough for our purposes here, but we would
like to emphasize that this description is only valid on the
subspace $\Arcn$ and cannot be generalized to the whole of $\Arc$.

\subsubsection{The dual graph} Informally the dual graph of an
element in $\Arcn$ is given as follows. The vertices are the
complementary regions.  Two vertices are joined by an edge if the
complementary regions border the same arc. Due to the orientation
of the surface this graph is actually a ribbon graph via the
induced cyclic order. Moreover the marked points on the boundary
make this graph into a marked ribbon graph. A more precise formal
definition is given in the next few paragraphs.

\subsubsection{Polygons and $\Arcn^0$}
By definition, in $\Arc^0_{\#}$ the complementary regions are
actually $k$-gons. Let $Poly(F,\Gamma,[i])$ be the set of these
polygons and let $Sides(F,\Gamma,[i])$ be the disjoint union of
sets of sides of the polygons. We define
$\del_{poly}:Sides(F,\Gamma,[i])\rightarrow Poly(F,\Gamma,[i])$ to
be the map which associates to a side $s$ of a polygon $p$ the
polygon $p$. The sides are either given by arcs or the boundaries.
We define the map $\lab:Sides(F,\Gamma,[i]) \rightarrow E_{\Gamma}
\bigcup V_{\Gamma}$ that associates the appropriate label. Notice
that for $e\in E_{\Gamma};|\lab^{-1}(e)|=2$ and for $v \in
V_{\Gamma}:|\lab^{-1}(v)|=1$.
 Thus there is a fixed point free involution
$\imath_{side}$ on the set $\lab^{-1}(E_{\Gamma})$ of sides of the
polygons marked by arcs which maps one side to the unique second
side carrying the same label. This in turn defines an involution
$\imath$ of pairs $(p,s)$ of a polygon together with a side in
$\lab^{-1} (E_{\Gamma})$ by mapping $s$ to $\imath_{side}(s)$ and
taking the polygon $p$ to the polygon $p':=\del_{poly}(\imath(s))$
of which $\imath_{side}(s)$ is a side. Although $p$ and $p'$ might
coincide the sides will differ making the involution $\imath$
fixed point free.

\subsubsection{Definition} For an
element $(F,\Gamma,[i],w) \in \Arc_{\#g}(n)$ we
define the dual graph
 to be the marked ribbon graph with a projective metric
 $(\hat \Gamma,ord,\hat w,\mk)$
 which is defined
 as follows. The vertices of $\hat\Gamma$ are the complementary
regions of the arc graph (i.e.\ the polygons) and the flags are
the pairs $(p,s)$ of a polygon (vertex) together with a side of
this polygon marked by an arc ($s\in \lab^{-1}(E_{\Gamma})$). The
map $\del$ is defined by $\del((p,s))=p$ and the involution
$\imath((p,s)):=(\del_{poly}(\imath_{side}(s)),\imath_{side}(s))$.

Each polygonal complimentary region is oriented by the orientation
induced by the surface, so that the sides of each polygon and thus
the flags of $\hat \Gamma$ at a given vertex $p$ have a natural
induced cyclic order $ord$ making $\hat \Gamma$ into a ribbon
graph.

 Notice that there is a one-one correspondence between edges of
the dual graph and edges of $\Gamma$. This is given by associating
to each edge $\{(p,s),\imath(p,s)\}$ the edge corresponding to the
arc $\lab(s)$.

We define a projective metric $\hat w$ on this graph by
associating to each edge $\{(p,s),\imath(p,s)\}$ the weight of the
arc labelling the side $s$ where $w$ is the projective metric on
the arc graph $\hat w(\{(p,s),\imath(p,s)\}):=w(\lab(s))$.

To define the marking, notice that the cycles of $\hat \Gamma$
correspond to the boundary components of the surface $F$. Let
$c_k$ be the cycle of the boundary component labelled by $k$. The
$k$-th boundary component lies in a unique polygon
$p=\del_{poly}(\lab^{-1}(k))$. Let $\prec_p$ be the cyclic order
on the set of sides of p, $\del^{-1}(p)$. Let $s_k$ be the side
corresponding to the boundary and let $N(s_k)$ the element
following $s_k$ in $\prec_p$. We define $\mk(c_k):=(p,N(s_k))$.

\subsubsection{Remark}
The above map will suffice for the purposes of this paper. For the
general theory and the reader acquainted with the constructions of
\cite{KLP}, the following Proposition will be helpful.

We do not, however, wish to go into technical details here on how
a marked weighted ribbon graph defines a configuration in the
sense of \cite{KLP} and refer the reader to \cite{hoch} for details.

\subsubsection{Proposition}\cite{hoch}\qua
{\sl For elements in $\Arcn$ the dual graph realizes the map $\Loop$.}

\subsection{The suboperads of $\Arc$ defined by the cacti operads}
\label{arcsection}

In this section, we would like to recall that the cacti without
spines and cacti can be embedded into the $\Arc$ operad up to an
overall scaling factor as defined in \cite{KLP}.
Moreover there is
an $S^1$ action given by the twist operator $\delta$. For the
complete details, we refer to \cite{KLP}.

In one direction the map is given by the dual graph discussed
above. In the other direction, the embedding is basically
constructed as follows: start by decomposing the cactus into the
arcs of its perimeter, where the break point of a cactus with
spines are the intersection points, the global zero and the local
zeros. Then one runs an arc from each arc to an outside pointed
circle which is to be drawn around the cactus configuration. The
arcs should be embedded starting in a counterclockwise fashion
around the perimeter of the circle. The marked points on the
inside circles which are the lobes of the cactus are the local
zeros for the cactus with spines and the global zero and the first
intersection point for a cactus without spines.

An equivalent formal definition in terms of graphs is given below.

\begin{figure}[ht!]
\epsfxsize = \textwidth \epsfbox{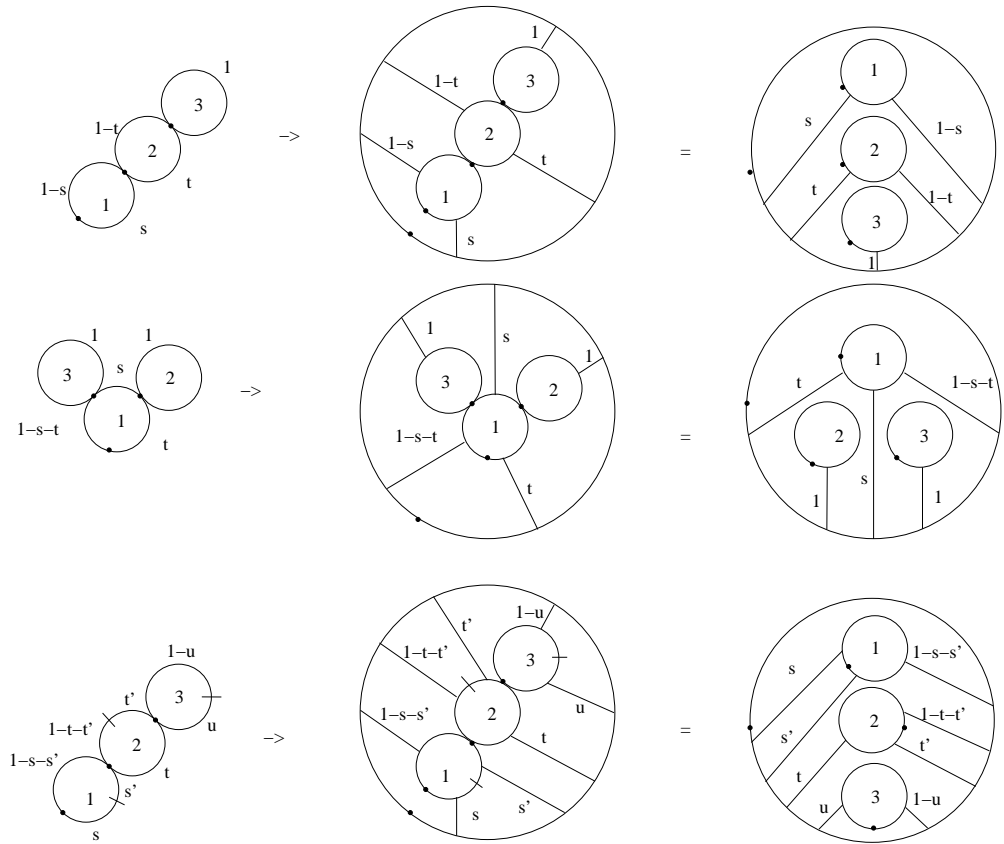}
\caption{\label{threeops}The embedding of cacti into arc}
\end{figure}

\begin{figure}[ht!]
\epsfxsize = \textwidth \epsfbox{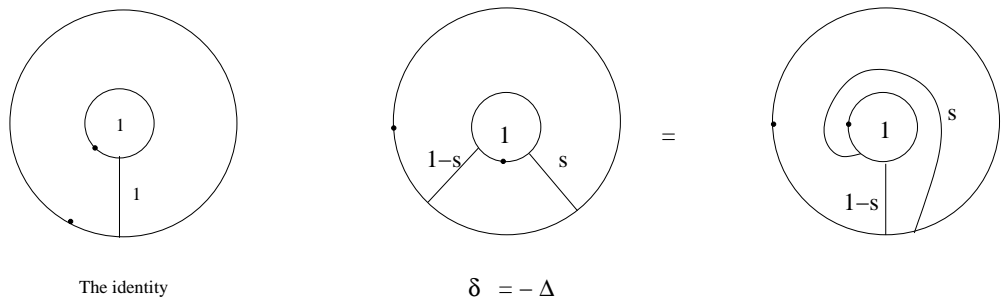}
\caption{\label{1ops}The identity and the twist $\delta$ yielding the BV
operator}
\end{figure}

For more orientation, we include two figures: the Figure
\ref{threeops} shows the framing i.e.\ embedding of two cacti
without spines and a cactus with spines into arc; Figure
\ref{1ops} shows the identity in arc and the family of weighted
arcs corresponding to the twist which yields the BV operator. The
Figures \ref{nospinesfig} V and \ref{spinesfig} V depict more
elaborate examples.

\subsubsection{Arcs and cacti}
The main result of the arc picture is summed up in the following Theorem.

\subsubsection{Theorem}{\sl There is an map of  (spineless) cacti
into $\Darc$ which maps $\Cacti$ bijectively onto
$\mathcal{D}\Tree_{cp}$ and $\Cact$ bijectively onto
$\mathcal{D}\Lintree_{cp}$. When restricted to its image this map
is an equivalence of operads.

Furthermore the suboperad in $\Darc$ generated by the
Fenchel-Nielsen type twist $\delta$ and the image of  spineless
cacti is equal to the image the cacti operad.

There are operadic  maps of the spineless cacti and the cacti
operad into $\Arc$ which are  equivalences when restricted to
their image.}

 \begin{proof} [Partial proof]
We will prove this claim on the level of sets and refer to
\cite{KLP} for the operad structure.

The map in one direction is given by associating the dual graph.
For an element in $\mathcal{D}\Tree_{cp}$ this graph is a marked
treelike ribbon graph with a metric, viz.\ a cactus. If the linear
orders agree, it is spineless.

For the map in the reverse direction, we realize a cactus $c$ as a
marked treelike ribbon graph with a metric $\Gamma(c)$. The
surface will be $\Sigma(\Gamma(c))$. The boundaries correspond to
the lobes and the outside circle and are hence labelled from $0$
to $n$.  Let $m(e)$ be the midpoint of the edge $e$. Then we
consider the arcs $m(e)\times [-1,1]$ on $\Sigma$. Notice that
these arcs come in a linear order at each boundary component
according to their linear order on the cactus and are labelled by
$w(e)\in \mathbb{R}_{>0}$.  Finally we mark off an interval on
each boundary such that the arcs on a boundary all end on this
interval and appear in the above linear order on the interval,
where the orientation of the interval is induced by that of the
surface. Contracting the interval to a point, we obtain the
desired element in $\Darc$.

The claim about $\Arc$ follows immediately from the equivalence of
operads $\Darc\rightarrow \Darc/\mathbb{R}_{>0}=\Arc$ which
contracts the factor $\mathbb{R}_{>0}$ of $\Darc$.
 \end{proof}

\subsubsection{Remark}
Alternatively, we can mark a point on each boundary such that the
arcs appear in their order on the complement of this point. This
alternative corresponds to the map called framing in \cite{KLP}
which we also used in our depictions in this paper.

\Addresses\recd

\end{document}